\documentstyle[amscd,amssymb,verbatim,12pt]{amsart}

\theoremstyle{plain} 
\newtheorem{thm}{Theorem}
\newtheorem{lem}[thm]{Lemma} 
\newtheorem{cor}[thm]{Corollary}
\newtheorem{prop}[thm]{Proposition}

\newenvironment{pfm}{\noindent{\em Proof of Theorem 
\ref{sets}.}}{\qed \\}

\newcommand{\ep}{\epsilon} 
\newcommand{\al}{\alpha}

\newcommand{\om}{\omega} 

\newcommand{\omoo}{{(\om^{\beta+1}+1)}}
\newcommand{\omo}{{(\om^{\beta+1})}}

\newcommand{\la}{\langle}
\newcommand{\ra}{\rangle}
\newcommand{\bs}{\backslash} 

\newcommand{\hr}{\hookrightarrow} 
\newcommand{\ol}{\overline}

\newcommand{\N}{{\Bbb N}}
\newcommand{\R}{{\Bbb R}}

\newcommand{\A}{\mbox{$\cal A$}}
\newcommand{\F}{\mbox{$\cal F$}}
\newcommand{\cS}{\mbox{$\cal S$}} 
\newcommand{\cP}{\mbox{$\cal P$}}
\newcommand{\Ps}{\mbox{${\cal P}_{<\infty}$}}
\newcommand{\Pb}{\mbox{${\cal P}_\infty$}}
\newcommand{\spn}{\operatorname{span}}
\newcommand{\supp}{\operatorname{supp}}
\newcommand{\sump}{\sum^p_{i=1}}

\begin{document}


\title[Compact subsets of $\Ps(\N)$]{Compact subsets of 
$\Ps(\N)$ with\\
applications to the embedding of\\
symmetric sequence spaces into $C(\alpha$)}
\author{Denny H.\ Leung} 
\address{Department of Mathematics\\ National
University of Singapore\\ Singapore 119260}
\email{matlhh@@math.nus.sg} 





\maketitle 

\begin{abstract} 
Let $\Ps(\N)$ be the set of all finite subsets of $\N$, endowed with 
the
product topology. A description of the compact subsets of $\Ps(\N)$ is 
given. Two applications of this result to Banach space theory are 
shown : (1) a characterization of the symmetric sequence spaces which 
embed into $C(\om^\om)$, and (2) a characterization, in terms of the 
Orlicz function $M$, of the Orlicz sequence spaces $h_M$ which embed 
into $C(K)$ for some countable compact Hausdorff space $K$. 
\end{abstract}


If $A$ is an arbitrary set, denote its power set by $\cP(A)$. 
Identifying $\cP(A)$ with $2^A$, and endowing it with the product 
topology, yields a compact Hausdorff topological space.  The symbols 
$\Ps(A)$ and $\Pb(A)$ stand for the subspaces consisting of all 
finite, respectively, all infinite subsets of $A$. In the first part 
of this paper, we study the compact subsets of $\Ps(\N)$. The fruit of 
this study is applied in the latter part to obtain some results 
concerning the embedding of symmetric Banach sequence spaces into 
$C(K)$, where $K$ is a countable compact Hausdorff space.

The main result in \S \ref{des} is Theorem \ref{sets}, which gives a 
description of the compact subsets of $\Ps(\N)$ in terms of a certain 
hierarchy of subsets $(\A^f_\beta)$ of $\Ps(\N)$ (see the definition 
in \S \ref{des} below). The motivation for the family $(\A^f_\beta)$ 
comes from the collection of ``admissible sets'' used in the 
definition of a classical counterexample (the Schrier space \cite{S}) 
in Banach space theory. Indeed, if $f$ is the identity function on 
$\N$, then $\A^f_1$ is precisely the collection of ``admissible sets'' 
used to define the Schrier space. Finite iterations of the 
construction appears in \cite{O}. It takes a certain amount of care, 
however, to extend the definition to transfinite ordinals. Careful 
choice is needed so that we may obtain the monotonicity condition in 
Theorem \ref{order}. The technical arguments, involving the behavior 
of ordinal numbers, are grouped together in \S \ref{bf}.

The remaining sections consist of certain applications of Theorem 
\ref{sets} to Banach space theory. Two main results are proved. The 
first, a combination of Theorems \ref{general} and \ref{converse}, is 
a  characterization of the symmetric Banach sequence spaces which 
embed into $C(\om^\om)$. In particular, it is observed that every 
Marcinkiewicz sequence space is isomorphic to a subspace of 
$C(\om^\om)$ (Corollary \ref{marc}). The other main result, Theorem 
\ref{last}, 
is the characterization, in terms of the Orlicz function $M$, of the 
Orlicz sequence spaces $h_M$ which embed into $C(\al)$ for some 
countable ordinal $\al$. The argument there involves the building of 
blocks (see the definition in \S \ref{subsp}) which are ``long'' with 
respect to the sets in $\A^f_\beta$.

For terms and notation concerning ordinal numbers and general 
topology, we refer to \cite{D}. The first infinite ordinal, 
respectively, the first uncountable ordinal, is denoted by $\om$, 
respectively, $\om_1$. Any ordinal is either $0$, a successor, or a 
limit. If $\al$ is a successor ordinal, denote its immediate 
predecessor by $\al - 1$. Following common practice, an ordinal $\al$ 
is identified with the set $\{\beta : \beta < \al\}$ of all of its 
predecessors. An ordinal is a Hausdorff topological space when endowed 
with the order topology. It is compact if and only if it is not a 
limit 
ordinal. If $K$ is a compact Hausdorff space, $C(K)$ denotes the space 
of all continuous real-valued functions on $K$. It is a Banach space 
under the norm $\|f\| = \sup_{t\in K}|f(t)|$. In particular, if $\al$ 
is an ordinal, we write $C(\al)$ for the space of all continuous 
real-valued functions on the compact Hausdorff space $\al + 1$. (The 
notation is inconsistent, but commonly used.) Detailed discussions of 
such spaces can be found in \cite{Se}. It is worth pointing out that 
every countable compact Hausdorff space is homeomorphic to some 
countable compact ordinal \cite{MS}. Thus every $C(K)$, where $K$ is 
countable compact, is isometric to some $C(\al)$, where $\al < \om_1$.
If $K$ is a topological space, its {\em derived set}\/ $K^{(1)}$ is 
the set of all of its limit points. A tranfinite sequence of derived 
sets may be defined as follows. Let $K^{(0)} = K$. If $\al$ is an 
ordinal, let $K^{(\al + 1)} = (K^{(\al)})^{(1)}$. Finally, for a limit 
ordinal $\al$, we define $K^{(\al)} = \cap_{\beta<\al}K^{(\beta)}$. 
The cardinality of a set $A$ is denoted by $|A|$. If $A$ and $B$ are 
nonempty subsets of $\N$, we say that $A < B$ if $\max A < \min B$. 
We also allow that $\emptyset < A$ and $A < \emptyset$ for any $A 
\subseteq \N$.

We follow standard Banach space terminology, as may be found in the 
book \cite{LT}. We say that a Banach space is a {\em sequence space}\/ 
if it is a vector subspace if the space of all real sequences. Such is 
the case, for instance, when a Banach space $E$ has a {\em (Schauder) 
basis}\/ $(e_k)$, i.e., every element $x \in E$ has a unique 
representation $x = \sum a_ke_k$ for some sequence of scalars $(a_k)$. 
Naturally, we identify every $x \in E$ with the sequence $(a_k)$  used 
in its representation. If $(e_k)$ is a basis of a Banach space $E$, 
there is a unique sequence of bounded linear functionals $(e'_k)$ on 
$E$ such that $\la e_j,e'_k\ra = 1$ if $j = k$, and $0$ otherwise. The 
sequence $(e'_k)$ is called the sequence of {\em biorthogonal 
functionals}\/ to the sequence $(e_k)$. It is a well known fact that 
every $x' \in E'$, the dual space of $E$, has a unique representation 
$x' = \sum a_ke'_k$, where the sum converges in the weak$^*$ topology 
on $E'$. Therefore, $E'$ may also be regarded as a sequence space. 
If $(e'_k)$ is a basis of $E'$ (so that the foregoing sum actually 
converges in norm for every $x' \in E'$), then the basis $(e_k)$ is 
said to be {\em shrinking}. If $x$ is an element of a sequence space, 
let $\supp x$ be the set of 
all coordinates $k$ at which $x$ is nonzero. Then we write $x < y$ to 
mean that $\supp x < 
\supp y$. The vector space consisting of all finitely supported real 
sequences is denoted by $c_{00}$. A basis $(e_k)$ of a Banach space is 
{\em 
unconditional}\/ if $\sum \ep_ka_ke_k$ converges for every choice of 
signs $(\ep_k)$  whenever $\sum a_ke_k$ converges. It is {\em 
symmetric}\/ if $\sum a_ke_{\pi(k)}$ converges for every permutation 
$\pi$ on $\N$ whenever $\sum a_ke_k$ converges. A symmetric basis is 
necessarily unconditional \cite[\S 3a]{LT}. We say that it is 
$1$-symmetric 
if $\|\sum \ep_ka_ke_{\pi(k)}\| = \|\sum a_ke_k\|$ for every choice of 
signs $(\ep_k)$, and every permutation $\pi$ on $\N$. Examples of 
Banach spaces with $1$-symmetric bases are $\ell^p$ $(1 \leq p < 
\infty)$, and $c_0$. These norms are defined by
\[
\|(a_k)\|_p = (\sum|a_k|^p)^{\frac{1}{p}} \\
\quad \text{and} \quad
\|(a_k)\|_\infty = \sup|a_k|
\]
respectively. 
A sequence $(x_k)$ in a Banach space is {\em normalized}\/ if $\|x_k\| 
= 1$ for all $k$. Given two sequences $(x_k)$ and $(y_k)$ in possibly 
different Banach spaces, we say that they are {\em equivalent} if 
there is a finite positive constant $C$ such that 
\[ C^{-1}\|\sum a_kx_k\| \leq \|\sum a_ky_k\| \leq C\|\sum a_kx_k\| \]
for every finitely supported sequence $(a_k)$.
Two Banach spaces $E$ and $F$ are said to be {\em isomorphic}\/ if 
they are linearly homeomorphic. We say that $E$ {embeds}\/ into $F$, 
$E \hr F$, if $E$ is isomorphic to a subspace of $F$.

\section{A description of the compact subsets of $\Ps(\N)$} 
\label{des}

Let $I$ be the collection of all countable limit ordinals. If $\al \in 
I$, denote by $I_\al$ the set $\{\beta \in I : \beta < \al\}$. 
Throughout this section, fix a limit ordinal
$\alpha < \om_1$. It is shown in \S \ref{bf} that there is a function 
$b_\al : 
I_\al\times\N \to \om_1$ such that 
\begin{enumerate}
\item for all $\gamma \in I_\al$, $(b_\al(\gamma,n))$ strictly 
increases to $\gamma$,
\item if $\gamma, \beta \in I_\al$, $n \in \N$, and $b_\al(\beta,n) < 
\gamma < \beta$, then $b_\al(\beta,n) < b_\al(\gamma,1)$.
\end{enumerate}
Let $\F$ be the collection of all functions $f :  \N \to \N$ which 
strictly increases to $\infty$. For each $f \in \F$, 
define subsets $\A^f_\beta$ of $\Ps(\N)$ inductively for all $\beta < 
\al$ as 
follows. Let $\A^f_0 = \{A \subseteq \N : |A| \leq 1\}$. If $\beta < 
\al$, 
and $\A^f_\gamma$ is defined for all $\gamma < \beta$, let 
\[ \A^f_{\beta} = \{A = \cup^k_{i=1}A_i : A_1 < \dots < A_k, \quad
A_1, \dots, 
A_k \in \A^f_{\beta-1}, \quad k \leq f(\min A)\} \]
if $\beta$ is a successor ordinal.
If $\beta$ is a limit ordinal, let
\[ \A^f_\beta = \{A : A \in \A^f_{b_\alpha(\beta,f(\min A))} \} . \]
It follows from the definition of the sets $\A^f_\beta$ and the 
properties of the function $b_\al$ that 
\begin{equation} \label{inclu}
\A^f_{b_\al(\beta,n)} \subseteq \A^f_{b_\al(\beta,n)+1} \subseteq 
\A^f_{b_\al(\beta,n+1)} 
\end{equation}
for all $\beta \in I_\al$, and all $n \in \N$. (See Theorem 
\ref{order} in \S \ref{bf}.)
We will use the sets $\A^f_\beta$ to describe the compact subsets of 
$\Ps(\N)$ up to the level of ``complexity'' $\al$. A subset $S$ of 
$\cP(\N)$ is said to be {\em hereditary}\/ if $A \subseteq B \in S$ 
implies that $A \in S$. We leave it to the reader to check that if $K$ 
is a hereditary subset of $\Ps(\N)$, then so is $K^{(\al)}$ for any 
ordinal $\al$.

\begin{prop} \label{her}
For each $\beta < \al$, $\A^f_\beta$ is a hereditary subset of 
$\Ps(\N)$. 
\end{prop}

\begin{pf}
Induct on $\beta$. The result clearly holds for $\beta = 0$. Now 
assume that $0 < \beta < \al$, and the result holds for all $\gamma < 
\beta$. Suppose that $B \subseteq A \in \A^f_\beta$. If $\beta$ is a 
successor, express $A$ as $\cup^k_{i=1}A_i$, where $A_1 < \dots < 
A_k$, $A_1, \dots, 
A_k \in \A^f_{\beta-1}$, and $k \leq f(\min A)$. Let $B_i = A_i \cap 
B$ for $1 \leq i \leq k$. Then $B_i \in \A^f_{\beta-1}$ since 
$\A^f_{\beta-1}$ is hereditary. Also, $f(\min A) \leq f(\min B)$. It 
follows easily that $B = \cup^k_{i=1}B_i \in \A^f_\beta$. If $\beta$ 
is a limit ordinal, 
then $A \in \A^f_{b_\al(\beta,f(\min A))}$. Since $b_\al(\beta,f(\min 
A)) < \beta$, $\A^f_{b_\al(\beta,f(\min A))}$ is hereditary. 
Therefore, 
$B \in A^f_{b_\al(\beta,f(\min A))}$. 
However, since $f(\min A) \leq f(\min B)$, we deduce from the relation 
(\ref{inclu}) that $A^f_{b_\al(\beta,f(\min A))} \subseteq 
A^f_{b_\al(\beta,f(\min B))}$. Thus $B \in A^f_{b_\al(\beta,f(\min 
B))}$, i.e., $B \in \A^f_\beta$. This completes the proof of the 
proposition.
\end{pf}

\begin{prop} \label{fg}
Let $f, g$ be two functions in \F, and suppose $k_0 \in \N$, $\beta < 
\al$. If $f(k) \geq g(k)$ for all $k \ge k_0$, and $A$ is an element 
of $\A^g_\beta$ such that $\min A \geq k_0$, then $A \in \A^f_\beta$.
\end{prop}

\begin{pf}
The proof is once again by induction on $\beta$. Since $\A^g_0 = 
\A^f_0$, the result holds for $\beta = 0$. Now suppose $0 < \beta < 
\al$, and the Proposition holds for all $\gamma < \beta$. Let $A \in 
\A^g_{\beta}$, $\min A \geq k_0$. If $\beta$ is a successor, then $A = 
\cup^k_{i=1}A_i$, where $A_1 < \dots < A_k$, 
$A_1, \dots, 
A_k \in \A^g_{\beta-1}$, and $k \leq g(\min A)$. Since $\min A_i \geq 
\min A \geq k_0$ for all $i$, $A_i \in \A^f_{\beta-1}$ by the 
inductive hypothesis. Also, observe that $f(\min A) \geq g(\min A)$.  
Hence the above decomposition of $A$ shows that $A \in \A^f_\beta$. On 
the other hand, if $\beta$ is a limit ordinal, then $A \in 
\A^g_{b_\al(\beta,g(\min A))}$. It follows from the relation 
(\ref{inclu}) and the fact that  $f(\min A) \geq g(\min A)$
that $\A^g_{b_\al(\beta,g(\min A))} \subseteq \A^g_{b_\al(\beta,f(\min 
A))}$. Therefore, $A \in \A^g_{b_\al(\beta,f(\min A))}$. As 
$b_\al(\beta,f(\min A)) < \beta$, we obtain from the inductive 
hypothesis that $A \in \A^f_{b_\al(\beta,f(\min A))}$. Hence $A \in 
\A^f_\beta$, as desired.
\end{pf}

\begin{thm} \label{sets}
Let $K$ be a compact, hereditary subset of $\Ps(\N)$ such that 
$K^{(\om^\beta+1)} = \emptyset$ for some ordinal $\beta < \al$. For 
all 
$C 
\in \Pb(\N)$, there exist $B \in \Pb(C)$, and a function $f \in \F$, 
such that $A \cap B \in \A^f_\beta$ for all 
$A \in 
K$.
\end{thm}

We prove Theorem \ref{sets} by induction on the 
ordinal $\beta$.  The main step, going from $\beta$ to $\beta +1$, 
requires
an induction of its own.  First, we introduce some more notation.
Recall the sets $\A^f_\beta$ defined above.  For any $m \in \N \cup
\{0\}$, let 
\[ \A^f_{\beta,m} = \{A = \cup^{k}_{i=1}A_i : A_1 < \dots < A_{k}, 
\quad A_1, \dots, A_{k} \in \A^f_\beta, \quad k \leq 3^m\} . \]
If $K \subseteq \Ps(\N)$, and $A \subseteq \N$, let $K_A = \{B \in K : 
A \subseteq B\}$. Also, for any subset $A$ of $\N$ and any $n \in \N$, 
let $A(<n)$, respectively $A(n)$, be the intersection of $A$ with the 
integer interval $[1,n)$, respectively, $[1,n]$.

\begin{lem} \label{tilde}
Let $K$ be a compact, hereditary subset of $\Ps(\N)$, and let $\gamma$ 
be a countable ordinal. If $A \in K \backslash K^{(\gamma)}$, and 
$\emptyset \in K^{(\gamma)}$, then there exists  
$a \in A$ such that $A(<a) \in K^{(\gamma)}$, and $A(a) \notin 
K^{(\gamma)}$. Moreover, $(K_A)^{(\gamma)} = \emptyset$.
\end{lem}

\begin{pf}
Write $A = \{n_1,\dots,n_l\}$, where $n_1 < n_2 < \dots < n_l$. Define 
$F_0 = \emptyset$, and $F_k = \{n_1,\dots,n_k\}$, $1 \leq k \leq l$. 
Let $k_0$ be the smallest integer such that $F_{k_0} \notin 
K^{(\gamma)}$. Such a $k_0$ exists since $F_{l} = A \notin 
K^{(\gamma)}$. Set $a = n_{k_0}$. Clearly, $A(<a) \in K^{(\gamma)}$, 
and $A(a) \notin K^{(\gamma)}$. Now if $B \in (K_A)^{(\gamma)}$, then 
$A \subseteq B 
\in K^{(\gamma)}$. Hence $A \in K^{(\gamma)}$, a 
contradiction.
\end{pf}

\begin{lem} \label{alpha} 
Suppose that Theorem \ref{sets} holds for some
ordinal $\beta < \al$.  If $K$ is a compact, hereditary subset of
$\Ps(\N)$ such that $K^{(\om^\beta\cdot 2^m+1)} = \emptyset$ for some 
$m \in \N \cup \{0\}$, then for all $C \in \Pb(\N)$, there exist $B 
\in \Pb(C)$,
and a function $f \in \F$ such that 
$A
\cap B \in \A^f_{\beta,m}$ for all $A \in K$.  
\end{lem}

\begin{pf}
We induct on $m$. When $m = 0$, the statement is simply Theorem 
\ref{sets} 
for the ordinal $\beta$. Now assume that the lemma holds for some $m 
\in 
\N \cup \{0\}$. Let $K$ be a compact, hereditary subset of $\Ps(\N)$ 
such that $K^{(\om^\beta\cdot 2^{m+1}+1)} = \emptyset$, and suppose 
that 
$C \in \Pb(\N)$ is given. We may obviously assume that 
$K^{(\om^\beta\cdot 2^{m}+1)} \neq \emptyset$. Since 
$(K^{(\om^\beta\cdot 
2^m)})^{(\om^\beta\cdot 2^m+1)} = \emptyset$, we may apply the 
inductive 
hypothesis to the set $K^{(\om^\beta\cdot 2^m)}$ to obtain a set $C_0 
\in 
\Pb(C)$, and a function $f_0 \in \F$ 
such that $A \cap C_0 \in \A^{f_0}_{\beta,m}$ whenever $A \in 
K^{(\om^\beta\cdot 2^m)}$. It remains to consider the behavior of the 
sets $A$ in $K \bs K^{(\om^\beta\cdot 2^m)}$.
Pick any $c_0 \in C_0$, and let $\cS_0 = \{S \in K \bs 
K^{(\om^\beta\cdot 2^m)} : \max S = c_0\}$. Clearly $\cS_0$ is finite. 
If $\cS_0 = \emptyset$, let $C_1 = C_0 \bs \{c_0\}$, and $f_1 = f_0$.
If $\cS_0 \neq \emptyset$, list its elements as $\cS_0 = \{S_1, \dots, 
S_l\}$. Using the fact that $(K_{S_i})^{(\om^\beta\cdot 2^m+1)} = 
\emptyset$, and the inductive hypothesis, we obtain sets 
\[ C_0 \bs \{c_0\} \supseteq C_0(S_1) \supseteq \dots \supseteq 
C_0(S_l)\]
in $\Pb(\N)$, and functions $g_1, \dots, g_l \in \F$ such that $A \cap 
C_0(S_i) \in \A^{g_i}_{\beta,m}$ whenever $A \in K_{S_i}$. Choose $f_1 
\in \F$ such that $f_1(k) \geq \max\{g_1(k), \dots, g_l(k)\}$ for all 
$k \in \N$, and let $C_1 = C_0(S_l)$. If $A \in K_{S_i}$ for some $i$, 
represent $A \cap C_0(S_i)$ as 
\[ A \cap C_0(S_i) = \cup^k_{j=1}D_j, \]
where $D_1 < \dots < D_{k}$, 
$D_1, \dots, D_{k} \in \A^{g_i}_\beta$, and $k \leq 3^m$. But since 
$f_1 \geq g_i$ for all $i$, $\A^{g_i}_\beta \subseteq \A^{f_1}_\beta$ 
by Proposition \ref{fg}. Also $\A^{f_1}_\beta$ is hereditary by 
Proposition \ref{her}. Hence $D_j \cap C_1 \in \A^{f_1}_\beta$ for all 
$j$. The representation $A \cap C_1 = \cup^k_{j=1}(D_j \cap C_1)$ 
shows that $A \cap C_1 \in \A^{f_1}_{\beta,m}$. Inductively, if $k 
\geq 1$, pick $c_k \in C_k$, $c_k > c_{k-1}$. If 
\[ \cS_k = \{S \in K \bs K^{(\om^\beta\cdot 2^m)} : \max S = c_k\} = 
\emptyset, \]
let $C_{k+1} = C_k \bs \{c_k\}$, and $f_{k+1} = f_k$.  Otherwise, 
repeating the argument above, we obtain $C_{k+1} \subseteq C_k \bs 
\{c_k\}$, and $f_{k+1} \in \F$, such that $A \in K_S$ for some $S \in 
\cS_k$ implies $A \cap C_{k+1} \in \A^{f_{k+1}}_{\beta,m}$.  Now let 
$B = \{c_0, c_1, c_2, \dots\}$. Then $B \in \Pb(C)$. Choose $f \in \F$ 
such  that $f(k) \geq \max\{f_0(k), f_1(k), \dots, f_k(k)\}$ for all 
$k$. Suppose $A \in K$. If $A \cap B \in K^{(\om^\beta\cdot 2^m)}$, 
then 
\[ A \cap B = A \cap B \cap C_0 \in \A^{f_0}_{\beta,m} \subseteq 
\A^{f}_{\beta,m} \]
by Proposition \ref{fg}. Since obviously $\A^{f}_{\beta,m} \subseteq 
\A^{f}_{\beta,m+1}$, the proof is complete in this case. If $A \cap B 
\notin K^{(\om^\beta\cdot 2^m)}$, apply Lemma \ref{tilde} to obtain 
$c_k \in A \cap B$ such that $(A\cap B)(<c_k) \in K^{(\om^\beta\cdot 
2^m)}$, but $(A\cap B)(c_k) \notin K^{(\om^\beta\cdot 2^m)}$. Let $D_1 
= (A\cap B)(<c_k)$, $D_2 = \{c_k\}$, and $D_3 = (A\cap B)\bs(D_1\cup 
D_2)$. Now 
$D_1 \in K^{(\om^\beta\cdot 2^m)}$ implies that $D_1 = D_1 \cap C_0 
\in \A^{f_0}_{\beta,m} \subseteq \A^{f}_{\beta,m}$. Also, $D_1 \cup 
D_2 \in 
\cS_k$, and $A \cap B \in K_{D_1\cup D_2}$. From the choices of 
$C_{k+1}$ and $f_{k+1}$, we conclude that 
\[ D_3 = A \cap B \cap C_{k+1} \in \A^{f_{k+1}}_{\beta,m} .\]
Since $\min D_3 \geq c_{k+1} \geq k+1$, and $f(j) \geq f_{k+1}(j)$ for 
all $j \geq k+1$, $D_3 \in \A^{f}_{\beta,m}$ by Proposition \ref{fg}.
Obviously, $D_2 \in \A^f_\beta$ as well. Since $A \cap B = D_1 \cup 
D_2 \cup D_3$, and $D_1 < D_2 < D_3$, it follows immediately that $A 
\cap B \in \A^f_{\beta,m+1}$. 
\end{pf}

The next result is the inductive step from $\beta$ to $\beta +1$ in 
Theorem \ref{sets}.

\begin{lem} \label{next}
If Theorem \ref{sets} holds for some ordinal $\beta < \al$, then 
it 
also holds for $\beta +1$.
\end{lem}

\begin{pf}
Suppose Theorem \ref{sets} holds for some ordinal $\beta < \al$. 
Then Lemma \ref{alpha} is applicable. Let $K$ be a compact, hereditary 
subset of $\Ps(\N)$ such that $K^\omoo = \emptyset$, and suppose that 
$C \in \Pb(\N)$ is given. If there exists $B \in \Pb(C)$ such that $B 
\cap (\cup_{A\in K}A) = \emptyset$, the result is trivial. We may thus 
assume otherwise. Since $K^{(\om^{\beta+1})}$ is finite, we can 
choose $c_1 \in C$, $c_1 > \max(\cup\{A : A \in K^\omo\})$, and $c_1 
\in \cup_{A\in K}A$. Now $(K_{\{c_1\}})^\omo = \emptyset$, hence 
$(K_{\{c_1\}})^{(\om^\beta\cdot 2^{m_1}+1)} = \emptyset$ for some $m_1 
\in \N$. By Lemma 
\ref{alpha}, there exist $C_1 \in \Pb(C\bs\{c_1\})$, and a function 
$f_1 \in \F$ such that $A \cap C_1 \in 
\A^{f_1}_{\beta,m_1}$ for all $A \in K_{\{c_1\}}$. If $C_n$ and $c_n$ 
have been chosen, let $c_{n+1} > c_n$, $c_{n+1} \in C_n$, $c_{n+1} \in 
\cup_{A\in K}A$. Then $c_{n+1} > c_1$; hence $(K_{\{c_{n+1}\}})^\omo = 
\emptyset$. Therefore $(K_{\{c_{n+1}\}})^{(\om^\beta\cdot 
2^{m_{n+1}}+1)} 
= \emptyset$ for some $m_{n+1} \in \N$. By Lemma \ref{alpha}, we 
obtain $C_{n+1} \in \Pb(C_n\bs\{c_{n+1}\})$, and a function $f_{n+1} 
\in \F$ such that $A \cap C_{n+1} \in 
\A^{f_{n+1}}_{\beta,m_{n+1}}$ for all $A \in K_{\{c_{n+1}\}}$. Let 
$B = \{c_1, c_2, \dots\}$, and let $f \in \F$ be such that $f(j) \geq 
\max\{f_1(j),\dots,f_j(j),2^{m_j}+1\}$ 
for all $j$. Then $B \in \Pb(C)$. Suppose $A \in K$. Let $\min(A\cap 
B) 
= c_j$. Then $A\cap B \in K_{\{c_j\}}$; hence $(A \cap B)\bs\{c_j\} =  
A \cap B \cap C_j \in \A^{f_j}_{\beta,m_j}$. Since $\min((A\cap 
B)\bs\{c_j\}) \geq c_{j+1} \geq j+1$, and $f(k) \geq f_j(k)$ for all 
$k \geq j+1$, $(A\cap B)\bs\{c_j\} \in \A^f_{\beta,m_j}$ by 
Proposition \ref{fg}. Therefore,
\[ (A \cap B) \bs \{c_j\} = \cup^k_{i=1}D_i, 
\quad D_1 < \dots < D_{k}, \quad D_1,\dots, D_{k} \in 
\A^f_\beta, \quad k \leq 2^{m_j} . \]
Let $E_1 = \{c_j\}$, $E_i = D_{i-1}$ for $2 \leq i \leq k+1$. 
Then 
\[ A \cap B = \cup^{k+1}_{i=1}E_i, \quad E_1 < \dots < 
E_{k+1}, 
\quad E_1, \dots, E_{k+1} \in \A^f_{\beta} . \]
Since $k+1 \leq 2^{m_j}+1 \leq f(j) \leq f(c_j) = f(\min(A\cap B))$, 
we 
conclude that $A \cap B \in \A^f_{\beta+1}$.
\end{pf}

\begin{pfm}
As indicated before, the proof is by induction on the 
ordinal $\beta$. First consider the case when $\beta = 0$. Let $K$ be 
a 
compact, hereditary subset of $\Ps(\N)$ such that $K^{(2)} = 
K^{(\om^0+1)} = \emptyset$. Also suppose that $C \in \Pb(\N)$ is 
given. We may assume that for all $B \in \Pb(C)$, $B \cap (\cup_{A\in 
K}A) \neq \emptyset$. Now $K^{(1)}$ is finite. Hence there exists $c_1 
\in C \cap (\cup_{A\in K}A)$ such that $c_1 > \max(\cup_{A\in 
K^{(1)}}A)$. For any $c \in C \cap (\cup_{A\in K}A)$, $c \geq c_1$, 
$(K_{\{c\}})^{(1)} = \emptyset$. Thus $K_{\{c\}}$ is finite. 
Therefore, 
there exists $c' \in C \cap (\cup_{A\in K}A)$ such that $c' > 
\max(\cup_{A\in K_{\{c\}}}A)$. Using this observation, we can choose a 
strictly increasing sequence $(c_n)$ in $C \cap (\cup_{A\in K}A)$  
such that $c_{n+1} > \max(\cup_{A\in K_{\{c_n\}}}A)$. Let $B = \{c_n : 
n \in \N\}$. Then $B \in \Pb(C)$. If $A \in K$, then $A \cap B \in K$. 
Let $\min (A\cap B) = c_j$. Since $A \cap B \in K_{\{c_j\}}$, and $c_n 
> \max(\cup_{A\in K_{\{c_j\}}}A)$ for all $n > j$, we see that $A \cap 
B = \{c_j\}$. Therefore, $A \cap B \in A^f_0$ for any $f \in \F$.

Now suppose $\beta < \al$, and the theorem has been 
proven for all ordinals $\gamma < \beta$. If $\beta$ is a successor 
ordinal, we 
appeal to Lemma \ref{next} to finish the proof. Otherwise, $\beta$ is 
a 
limit ordinal. Let $K$ be a compact, hereditary subset of $\Ps(\N)$ 
such that $K^{(\om^\beta+1)} = \emptyset$, and let $C \in \Pb(\N)$ be 
given. Once again, we assume without loss of generality that $B \cap 
(\cup_{A\in K}A) \neq \emptyset$ for all $B \in \Pb(C)$. Since 
$K^{(\om^\beta)}$ is finite, there exists $c_1 \in C \cap (\cup_{A\in 
K}A)$ such that $c_1 > \max(\cup\{A : A \in K^{(\om^\beta)}\})$. Then 
$(K_{\{c_1\}})^{(\om^\beta)} = \emptyset$. Hence 
$K^{(\om^{b_\al(\beta,n_1)}+1)} = \emptyset$ for some $n_1 \in \N$.
As $b_\al(\beta,n_1) < \beta$, we may apply the inductive hypothesis 
to 
obtain $C_1 \in \Pb(C\bs\{c_1\})$ and a function $f_1 \in \F$ such 
that $A \cap C_1 \in 
\A^{f_1}_{b_\al(\beta,n_1)}$ for all $A \in K_{\{c_1\}}$.  If $C_k$ 
and 
$c_k$ have been chosen, let $c_{k+1} \in C \cap (\cup_{A\in K}A)$ be 
such that $c_{k+1} > c_k$. Then $(K_{\{c_{k+1}\}})^{(\om^\beta)} = 
\emptyset$. Hence there exist $n_{k+1} \in \N$, $C_{k+1} \in 
\Pb(C_k\bs\{c_{k+1}\})$, and a function $f_{k+1} \in \F$ such that 
\[ A \cap C_{k+1} \in \A^{f_{k+1}}_{b_\al(\beta,n_{k+1})} \quad 
\text{for all} \quad A \in K_{\{c_{k+1}\}}. \]
Now let $B = \{c_1, c_2, \dots\}$, and pick $f \in \F$ so that 
$f(j) \geq 
\max\{f_1(1),\dots,f_j(j),n_j+1\}$ for all $j$. Suppose $A \in K$, and 
$\min(A\cap B) = c_k$. As $A \cap B \in K_{\{c_k\}}$, 
\[ (A \cap B)\bs \{c_k\} = A \cap B \cap C_k \in 
\A^{f_k}_{b_\al(\beta,n_k)}. \]
Since $\min((A\cap B)\bs\{c_k\}) \geq c_{k+1} \geq k+1$, and $f(j) 
\geq f_k(j)$ for all $j \geq k+1$, $(A\cap B)\bs\{c_j\} 
\in \A^f_{b_\al(\beta,n_k)}$ by Proposition \ref{fg}. Now let $D_1 = 
\{c_j\}$, $D_2 = (A\cap 
B) \bs \{c_j\}$. Then $D_1 < D_2$, and $D_1, D_2 \in 
\A^f_{b_\al(\beta,n_k)}$. Also, $f(\min(A\cap B)) \geq f(1) \geq 2$. 
Therefore, $A \cap B = D_1 \cup D_2 \in \A^f_{b_\al(\beta,n_k)+1}$. 
Finally, since $f(\min(A\cap B)) = f(c_k) \geq 
n_k +1$, if follows from the relation (\ref{inclu}) that
$\A^f_{b_\al(\beta,n_k)+1} \subseteq 
\A^f_{b_\al(\beta,f(\min(A\cap B)))}$. Hence $A \cap B \in 
\A^f_{b_\al(\beta,f(\min(A\cap B)))}$, i.e., $A \cap B \in 
\A^f_\beta$.
This completes the proof of the theorem. 
\end{pfm}

\noindent{\em Remark}. We leave it to the reader to check that 
$\A^f_\beta$ is a compact, hereditary subset of $\Ps(\N)$ such that 
$(\A^f_\beta)^{(\om^\beta+1)} = \emptyset$. Therefore, the description 
of $K$ given in Theorem \ref{sets} is sharp except for the fact that 
we have to pass to an infinite subset $B$.

\section{Construction of the function $b_\al$} \label{bf}

In this section, we give the construction of the function $b_\al$ used 
in \S \ref{des}. The crucial step is given in Lemma \ref{c}.

\begin{prop} \label{h}
There is a continuous order isomorphism $h$ from $\om_1$ onto $I$.
\end{prop}

\begin{pf}
Since $I$ is a cofinal subset of $\om_1$, there is an order 
isomorphism $h$ from $\om_1$ onto $I$. We claim that $h$ is 
continuous. Indeed, let $(\al_n)$ be a sequence of countable ordinals 
which strictly increases to some $\al$. Then $(h(\al_n))$ strictly 
increases to some $\beta < \om_1$. Clearly, $\beta \in I$. Say $\beta 
= h(\gamma)$. Since $h(\al) \geq h(\al_n)$ for all $n$, $h(\al) \geq 
h(\gamma)$. Hence $\al \geq \gamma$. But if $\gamma < \al$, then there 
exists $n$ such that $\gamma < \al_n$. This implies that $\beta = 
h(\gamma) < h(\al_n)$, which is impossible. Therefore, $\gamma = \al$. 
This shows that $h$ is continuous on $\om_1$.
\end{pf}

\begin{lem} \label{c}
For every $\al \in I$, there is a function $c_\al : I_\al \times \N 
\to \om_1$ such that for all $\gamma \in I_\al$,
\begin{enumerate}
\item the sequence $(c_\al(\gamma,n))$ strictly increases to $\gamma$,
\item $\sup\{c_\al(\beta,n) : (\beta,n) \in I_\al\times \N,   
c_\al(\beta,n) < \gamma < \beta\} < \gamma$.
\end{enumerate}
(Here, we let $\sup\emptyset = 0$.)
\end{lem}

\begin{pf}
We induct on the elements of $I$. The result holds vacuously for $\min 
I = \om$, since $I_\om = \emptyset$. Now suppose $\om < \al_0 \in I$, 
and functions $c_\al$ are found for 
all $\al \in I_{\al_0}$ satisfying conditions (1) and (2). We want to 
show that the result also holds 
for $\al_0$. Using the function $h$ from Proposition \ref{h}, write 
$\al_0 = h(\beta_0)$. Then $0 < \beta_0 < \om_1$. We consider 3 
separate cases.\\

\noindent\underline{Case 1}. $\beta_0$ is a successor and $\beta_0 - 
1$ is not a limit.\\
Define $c_{\al_0} : I_{\al_0} \times \N \to \om_1$ by 
\[ c_{\al_0}(\gamma,n) = \begin{cases}
   c_{h(\beta_0 - 1)}(\gamma,n)& \text{if $\gamma \in I_{h(\beta_0 - 
1)}$},\\
   h(\beta_0 - 2) + n & \text{if $\gamma = h(\beta_0 - 1)$}.
   \end{cases} \]
Here, we let $h(\beta_0 - 2) = 0$ if $\beta_0 - 1 = 0$.   
If $\gamma \in I_{h(\beta_0 - 1)}$, then condition (1) holds by the 
inductive choice of $c_{h(\beta_0 - 1)}$. Consider condition (2). If  
$\beta = h(\beta_0 - 1)$, then $c_{\al_0}(\beta,n) > \gamma$. 
Therefore
\begin{multline*}
\sup\{c_{\al_0}(\beta,n) : (\beta,n) \in I_{\al_0}\times \N,   
c_{\al_0}(\beta,n) < \gamma < \beta\}\\
= \sup\{c_{h(\beta_0 - 1)}(\beta,n) : (\beta,n) \in 
I_{h(\beta_0-1)}\times \N,
c_{h(\beta_0 - 1)}(\beta,n) < \gamma < \beta\} <  
\gamma
\end{multline*}
by the inductive choice of $c_{h(\beta_0 - 1)}$. Thus condition (2) 
holds if $\gamma \in I_{h(\beta_0 - 1)}$. Finally, if $\gamma = 
h(\beta_0 - 
1)$, then clearly $\gamma = h(\beta_0 - 2) + \om$. Thus condition (1) 
holds. Moreover, there is no $\beta \in I_{\al_0}$ such that $\beta > 
\gamma$. So condition (2) also holds.\\

\noindent\underline{Case 2}.  $\beta_0$ is a successor and $\beta_0 - 
1$ is a limit.\\
Choose a sequence of ordinals $(\tau_n)$ which strictly increases to 
$\beta_0 - 1$.
Define $c_{\al_0} : I_{\al_0} \times \N \to \om_1$ by 
\[ c_{\al_0}(\gamma,n) = \begin{cases}
   c_{h(\beta_0 - 1)}(\gamma,n)& \text{if $\gamma \in I_{h(\beta_0 - 
1)}$},\\
   h(\tau_n)& \text{if $\gamma = h(\beta_0 - 1)$}.
   \end{cases} \]
It is clear by the definition of $c_{\al_0}$ that condition (1) holds. 
If $\gamma \in I_{h(\beta_0 - 1)}$, 
\begin{equation*}
\begin{split}
\sup\{&c_{\al_0}(\beta,n) : (\beta,n) \in I_{\al_0}\times\N, 
c_{\al_0}(\beta,n) < \gamma < \beta\} = \\
&\max\bigl\{\sup\{c_{h(\beta_0 - 1)}(\beta,n) : (\beta,n) \in 
I_{h(\beta_0-1)}\times\N, c_{h(\beta_0 - 1)}(\beta,n) < 
\gamma < \beta\},\\
&\hspace{5.6ex}\sup\{h(\tau_n) : n \in \N, h(\tau_n) < 
\gamma\}\bigr\}. 
\end{split}
\end{equation*}
The first supremum on the right is $< \gamma$ by the inductive choice 
of $c_{h(\beta_0 - 1)}$. The second supremum on the right is also $< 
\gamma$ because $(h(\tau_n))$ strictly increases to $h(\beta_0 - 1)$, 
and $\gamma < h(\beta_0 - 1)$. Thus condition (2) is verified if 
$\gamma \in I_{h(\beta_0 - 1)}$. On the other hand, if $\gamma = 
h(\beta_0 - 1)$, there is no $\beta \in I_{\al_0}$ such that $\beta > 
\gamma$. Hence condition (2) is fulfilled trivially.\\

\noindent\underline{Case 3}. $\beta_0$ is a limit ordinal.\\
Choose a sequence of ordinals $(\tau_k)$ which increases to $\beta_0$, 
such that $\tau_k + 2 < \tau_{k+1}$ for all $k \in \N$. If $k \in \N$, 
$\gamma \in I$, and $h(\tau_k+2) \leq \gamma \leq h(\tau_{k+1})$, then 
$\gamma \in I_{h(\tau_{k+1}+1)}$. Since $h(\tau_{k+1} + 1) \in 
I_{\al_0}$, by the inductive hypothesis,
$c_{h(\tau_{k+1}+1)}(\gamma,n)$
increases to $\gamma$ as $n \to \infty$. But $\gamma \geq h(\tau_k+2) 
> h(\tau_k +1)$. Therefore, there exists $n_\gamma \in \N$ such that 
$c_{h(\tau_{k+1}+1)}(\gamma,n) > h(\tau_k +1)$ if $n > n_\gamma$. 
Define $c_{\al_0} : I_{\al_0} \times \N \to \om_1$ as 
follows: 
\[ c_{\al_0}(\gamma,n) = \begin{cases}
   c_{h(\tau_1 + 1)}(\gamma,n)& \text{if $\gamma \leq h(\tau_1)$},\\
   h(\tau_k)+n & \text{if $\gamma = h(\tau_k+1)$, $k \in \N$},\\
   c_{h(\tau_{k+1}+1)}(\gamma,n+n_\gamma)& \text{if $h(\tau_k+2) \leq  
      \gamma \leq h(\tau_{k+1})$, $k \in \N$}.
   \end{cases} \]
Let $\gamma \in I_{\al_0}$, if  $\gamma \leq h(\tau_1)$, or if 
$h(\tau_k+2) \leq \gamma \leq h(\tau_{k+1})$ for some $k \in \N$. Then 
$(c_{\al_0}(\gamma,n))$ strictly increases to $\gamma$ by the 
inductive choices of $c_{h(\tau_1+1)}$ and $c_{h(\tau_{k+1}+1)}$ 
respectively. On the other hand, if $\gamma = h(\tau_k+1)$ for some $k 
\in \N$, then $(c_{\al_0}(\gamma,n))$ increases to $\gamma$, since 
$\gamma = h(\tau_k+1) = h(\tau_k) + \om$. This verifies condition (1) 
for the function $c_{\al_0}$. 
To check condition (2), first observe that if $\beta \in I_{\al_0}$, 
and $h(\tau_k+2) \leq \beta \leq h(\tau_{k+1})$ for some $k \in \N$, 
then $c_{\al_0}(\beta,n) > h(\tau_k +1)$ by the choice of $n_\beta$. 
It is also clear that if $\beta = h(\tau_k+1)$ for some $k$, then 
$c_{\al_0}(\beta,n)> h(\tau_k)$. Using these observations, we deduce 
that if $\gamma, \beta \in I_{\al_0}$, and 
$c_{\al_0}(\beta,n) < \gamma < \beta$, then
\[ \beta \leq \begin{cases}
   h(\tau_1)& \text{if $\gamma \leq h(\tau_1)$},\\
   h(\tau_{k+1})& \text{if $h(\tau_k+2) \leq \gamma \leq               
    h(\tau_{k+1})$},
   \end{cases} \]
while no such $\beta$ exists if $\gamma = h(\tau_k+1)$ for some $k$.  
Therefore, given $\gamma \in I_{\al_0}$,
\begin{equation*}
\begin{split}
\{c_{\al_0}&(\beta,n) : (\beta,n) \in I_{\al_0}\times\N,  
c_{\al_0}(\beta,n) < \gamma < \beta\}\\ 
&= \begin{cases}
   \{c_{h(\tau_1+1)}(\beta,n): (\beta,n) \in I_{h(\tau_1+1)}\times\N,  
     c_{h(\tau_1+1)}(\beta,n) < \gamma < \beta\},& \\
   \emptyset, & \\
   \{c_{h(\tau_{k+1}+1)}(\beta,n+n_\beta): (\beta,n) \in               
     I_{h(\tau_{k+1}+1)}\times\N, c_{h(\tau_{k+1}+1)}(\beta,n+n_\beta) 
     < \gamma < \beta\}.&        
   \end{cases}       
\end{split}
\end{equation*}
according to whether $\gamma \leq h(\tau_1)$, $\gamma = h(\tau_k+1)$, 
or $h(\tau_k+2) \leq \gamma \leq h(\tau_{k+1})$ for some $k$. Using 
the inductive properties of the functions $c_{h(\tau_1+1)}$ and 
$c_{h(\tau_{k+1}+1)}$, we 
deduce easily that $c_{\al_0}$ satisfies condition (2).
\end{pf}

\begin{thm}\label{b}
For every $\al \in I$, there is a function $b_\al : I_\al \times \N 
\to \om_1$ such that 
\begin{enumerate}
\item for all $\gamma \in I_\al$, $(b_\al(\gamma,n))$ strictly 
increases to $\gamma$,
\item if $\gamma, \beta \in I_\al$, $n \in \N$, and $b_\al(\beta,n) < 
\gamma < \beta$, then $b_\al(\beta,n) < b_\al(\gamma,1)$.
\end{enumerate}
\end{thm}

\begin{pf}
Given $\al \in I$, obtain a function $c_\al : I_\al \times \N 
\to \om_1$ using Lemma \ref{c}. If $\gamma \in 
I_\al$, let 
\[ \tau_\gamma = \sup\{c_\al(\beta,n) : (\beta,n) \in I_\al\times\N, 
c_\al(\beta,n) < \gamma < \beta\}. \]
Then $\tau_\gamma < \gamma$ by Lemma \ref{c}. Also 
$(c_\al(\gamma,n))$ increases to $\gamma$. Hence,  there exists 
$m_\gamma \in \N$ such that $c_\al(\gamma,n) > \tau_\gamma$ if $n > 
m_\gamma$. Define $b_\al : I_\al \times \N 
\to \om_1$ by $b_\al(\gamma,n) = c_\al(\gamma, 
n+m_\gamma)$. Condition (1) of the theorem follows from condition (1) 
of Lemma \ref{c}. Now suppose $\gamma, \beta \in I_\al$, $n \in 
\N$, and $b_\al(\beta,n) < \gamma < \beta$. Then
\[ c_\al(\beta,n+m_\beta) = b_\al(\beta,n) < \gamma < \beta . \]
It follows from the definition of $\tau_\gamma$ that $b_\al(\beta,n) 
\leq \tau_\gamma$. But 
\[ \tau_\gamma < c_\al(\gamma,1+m_\gamma) = b_\al(\gamma,1). \]
This proves condition (2).
\end{pf}

From now on, fix a collection of functions $\{b_\al : \al \in I\}$ 
satisfying the conditions of Theorem \ref{b}.

\begin{cor} \label{p}
Let $\al \in I$.  Define $p : \al \backslash \{0\} \to \om_1$ by
\[ p(\beta) = \begin{cases}
     b_\al(\beta,1)& \text{if $\beta$ is a limit ordinal},\\
     \beta - 1& \text{if $\beta$ is a successor ordinal}.
     \end{cases}\]
If $\beta \in I_\al$, and $n \in \N$, there exists $k \in \N\cup\{0\}$ 
such that $b_\al(\beta,n) + 1 = p^k(b_\al(\beta,n+1))$.
\end{cor}     

\begin{pf}
Clearly, if $p^k(b_\al(\beta,n+1)) \neq 0$, then 
$p^k(b_\al(\beta,n+1)) > p^{k+1}(b_\al(\beta,n+1))$. Let 
\[ J = \{k \in \N\cup\{0\} : p^k(b_\al(\beta,n+1)) \geq b_\al(\beta,n) 
+ 1\}. \]
Obviously, $0 \in J$. Since every strictly decreasing sequence of 
ordinals is finite, $J$ is a finite set. Let $j_0 = \max J$, and 
denote $p^{j_0}(b_\al(\beta,n+1))$ by $\gamma$. By definition, $\gamma 
\geq b_\al(\beta,n) + 1$. If $\gamma$ is a successor, then $p(\gamma) 
= \gamma - 1 < b_\al(\beta, n) + 1$. Hence $\gamma = b_\al(\beta, n) + 
1$. On the other hand, if $\gamma$ is a limit ordinal, then$\ 
b_\al(\beta,n) < \gamma < \beta$. Therefore, $b_\al(\beta,n) < 
b_\al(\gamma,1) = p(\gamma)$ by the properties of $b_\al$. But then 
$j_0 + 1 \in J$, contrary to the choice of $j_0$. Thus, it must be 
that $\gamma = p^{j_0}(b_\al(\beta,n+1)) = b_\al(\beta,n) + 1$.
\end{pf}

\begin{thm} \label{order}
If $\al \in I$, $\beta \in I_\al$, and $n \in \N$, then for any $f \in 
\F$,
\[ \A^f_{b_\al(\beta,n)} \subseteq \A^f_{b_\al(\beta,n)+1} \subseteq 
\A^f_{b_\al(\beta,n+1)} .\]
\end{thm}

\begin{pf}
Observe that if $\gamma \in I_\al$, then $\A^f_\gamma \subseteq 
\A^f_{\gamma+1}$. We prove the second inclusion by induction on the 
elements of 
$I_\al$. First consider the case when $\beta = \min I_\al = \om$. Then 
every $b_\al(\beta,m)$ is a finite ordinal. Thus, there exists $j \in 
\N\cup\{0\}$ such that $b_\al(\beta,n+1) = b_\al(\beta,n)+1+j$. It 
follows readily from the above observation that 
$\A^f_{b_\al(\beta,n)+1} \subseteq \A^f_{b_\al(\beta,n+1)}$. Now 
suppose $\beta \in I_\al$, and $\A^f_{b_\al(\gamma,n)+1} \subseteq 
\A^f_{b_\al(\gamma,n+1)}$ for all $\gamma \in I_\beta$. 
We claim that if $\gamma$ is {\em any}\/ ordinal such that $0 < \gamma 
< \beta$, then $\A^f_{p(\gamma)} 
\subseteq \A^f_\gamma$, where $p$ is the function defined in Corollary 
\ref{p}.  Indeed, if $\gamma$ is a successor, the statement is simply 
the observation made at the beginning of the proof. On the other hand, 
if $\gamma$ is a limit ordinal, then
\[ \A^f_{p(\gamma)} = \A^f_{b_\al(\gamma,1)} \subseteq 
\A^f_{b_\al(\gamma,1)+1} \subseteq \A^f_{b_\al(\gamma,2)} 
\subseteq \A^f_{b_\al(\gamma,2)+1} \subseteq \A^f_{b_\al(\gamma,3)}
\subseteq \dotsm . \]
Hence $\A^f_{p(\gamma)} \subseteq \A^f_{b_\al(\gamma,m)}$ for all $m$. 
In particular, if $A \in \A^f_{p(\gamma)}$, then $A \in 
\A^f_{b_\al(\gamma,f(\min A))}$, i.e., $A \in \A^f_\gamma$. This 
proves the claim. By Corollary \ref{p}, there exists $k \in 
\N\cup\{0\}$ such that $b_\al(\beta,n)+1 = p^k(b_\al(\beta,n+1))$. 
Clearly, $p^j(b_\al(\beta,n+1)) < \beta$ whenever $0 \leq j \leq k$. 
Thus we can conclude from the claim that 
\[ \A^f_{b_\al(\beta,n)+1} = \A^f_{p^k(b_\al(\beta,n+1))} \subseteq 
\A^f_{b_\al(\beta,n+1)}, \]
as required.
\end{pf}

\section{An embedding result}

In this section, we prove a general result which shows that members of 
a certain class of symmetric sequence spaces can be embedded into 
$C(\om^\om)$. This class includes the Marcinkiewicz sequence spaces. 
We also prove that it includes all Orlicz sequence spaces $h_M$ such 
that $\lim_{t\to 0}M(\eta t)/M(t) = 0$ for some $\eta > 0$.
A norm $\rho$ defined on a vector lattice $E$ is a {\em lattice
norm}\/ if $\rho(x) \leq \rho(y)$ whenever $x, y \in E$ and $|x|
\leq |y|$.  We also say that a norm $\rho$ on $\R^n$ is {\em
normalized}\/ if $\|\cdot\|_\infty \leq \rho \leq \|\cdot\|_1$. If $a 
= (a_n)$ is a finite real sequence, or an infinite real sequence which 
converges to $0$, we let $a^* = (a^*_n)$ denote the decreasing 
rearrangement of the sequence $(|a_n|)$.

\begin{thm} \label{general} 
Suppose that $(t_n)$ is a strictly positive real sequence
which decreases to $0$, and for each $n \in \N$, $\rho_n$ is
a normalized lattice norm on $\R^n$.  If a Banach space $E$ has a 
basis $(e_k)$ such that 
\begin{equation} \label{norm} 
\|\sum
a_ke_k\| = \sup_n t_n\,\rho_n(a^*_1,\ldots,a^*_n) 
\end{equation} 
for
every element $\sum a_ke_k \in E$, then $E$ is isomorphic to a
subspace of $C(\om^\om)$.  
\end{thm}

\begin{pf} 
For each $n$, denote by $B_n$ the unit ball of the dual of
$(\R^n,\rho_n)$.  Since $\rho_n$ is normalized, $\|b\|_\infty \leq 1$
if $b \in B_n$.  Choose a finite subset $A_n$ of $B_n$ consisting of
decreasing sequences so that if $a \in \R^n$, then 
\begin{equation} \label{An} 
\frac{1}{2}\rho_n(a^*) \leq \sup_{b\in A_n}\la a^*,b\ra
\leq \rho_n(a^*) .  
\end{equation} 
Let $K_n$ be the set of all $c
\in \R^\N$ such that either $c = 0$, or there are a $b = (b_k) \in
A_n$, and $1 \leq k_1 < k_2 < \cdots < k_i \leq n$ such that $c^* =
t_n(b_{k_1},\ldots,b_{k_i},0,\ldots)$.  Endow $\R^\N$ with the
topology of pointwise convergence.  Then each $K_n$ is a compact
subset.  Let $\cP_n(\N) = \{A \subseteq \N : |A| \leq n\}$.
Since $(\cP_n(\N))^{(m)} = \emptyset$ for $m > n$,
and $A_n$ is finite, we see that $(K_n)^{(m)} = \emptyset$ if $m > n$.
Let $K = \cup^\infty_{n=1}K_n$.  We claim that $K$ is pointwise
compact.  Indeed, if $c \in
\ol{K}\backslash\cup^\infty_{n=1}\ol{K_n}$, then for any $n, k \in
\N$, and all $\ep > 0$, there exists $d = (d_k) \in
K\backslash\cup^n_{j=1}K_j$ such that $|c_k - d_k| < \ep$.  Hence 
\[
|c_k| < |d_k| + \ep \leq \|d\|_\infty + \ep \leq t_{n+1} + \ep .\] 
Since $k, n$ and $\ep$ are arbitrary, $c_k = 0$ for all $k \in \N$.
But then $c = 0 \in \cup^\infty_{n=1}K_n$, a contradiction.  Thus 
\[
\ol{K} = \cup^\infty_{n=1}\ol{K_n} = \cup^\infty_{n=1}K_n = K .  \] 
So $K$ is pointwise closed.  Since $K$ is coordinatewise bounded, it
is pointwise compact.  Repeating the above proof shows that 
\[
K^{(m)} = \cup^\infty_{n=1}(K_n)^{(m)} = \cup^\infty_{n=m}(K_n)^{(m)}
\subseteq \cup^\infty_{n=m}{K_n} .  \] 
Therefore, 
\[ K^{(\om)} =
\cap^\infty_{m=1}K^{(m)} \subseteq
\cap^\infty_{m=1}(\cup^\infty_{n=m}{K_n}) = \{0\} .  \] 
Denote by
$(e'_k)$ the sequence in $E'$ biorthogonal to $(e_k)$.  Identify each
element $a = (a_k)$ in $K$ with the functional $\sum a_ke'_k$ in $E'$.
Because of (\ref{norm}) and (\ref{An}), $K$ (thus identifed) is a 
subset of the unit ball of $E'$. Since $K$ is pointwise compact, it is 
a $\sigma(E',E)$-compact subset of $E'$.  Moreover, $K^{(\om)} 
\subseteq
\{0\}$ in the $\sigma(E',E)$-topology.  By a result of Semadeni (cf.\
\cite[Corollary 5.2]{La}), we see that $K$ is homeomorphic to $\alpha
+1$ for some ordinal $\alpha \leq \om^\om$.  It follows from 
(\ref{norm}) and (\ref{An}) that $E$ embeds
into $C(K)$.  Therefore, 
\[ E \hr C(\alpha) \hr C(\om^\om) , \] 
as
required.  
\end{pf}

In Theorem \ref{general}, if each $\rho_n$ is the $\ell^1$-norm, then 
the resulting space is called a {\em Mar\-cin\-kie\-wicz sequence 
space}. If, in addition, we let $t_n = n^{\frac{1}{p}-1}$, $1 < p < 
\infty$, then we obtain the closed linear span of the unit vectors in 
weak $\ell^p$.

\begin{cor} \label{marc}
Every Marcinkiewicz sequence space is isomorphic to a subspace of 
$C(\om^\om)$. In particular, the closed linear span of the unit 
vectors in weak $\ell^p$  is isomorphic to a subspace of $C(\om^\om)$.
\end{cor}

Next, we turn to the class of Orlicz sequence spaces. 
An Orlicz function is an increasing, convex function $M : [0,\infty) 
\to \R$ such that $M(0) = 0$ and $M(1) = 1$. It is {\em 
nondegenerate}\/ if $M(t) > 0$ for all $t > 0$. The corresponding {\em 
Orlicz sequence space}\/ $h_M$ is the space of all real sequences 
$(a_k)$ such that $\sum_kM(c|a_k|) < \infty$ for all $c < \infty$, 
endowed with the norm
\[ \|(a_k)\| = \inf\{\rho > 0 : \sum^\infty_{k=1}M(|a_k|/\rho) \leq 
1\} . \]
It is well known that the coordinate unit vectors form a normalized 
$1$-symmetric 
basis of $h_M$. Moreover, the dual $h'_M$ of $h_M$ can be identified 
with the collection of all real sequences $(b_k)$ such that $\la(a_k), 
(b_k)\ra = \sum a_kb_k < \infty$ for all $(a_k) \in h_M$, $\|(a_k)\| 
\leq 1$. For further results concerning Orlicz spaces, we refer to 
\cite{LT}.

We will show 
that Theorem \ref{general} is applicable to any Orlicz sequence space 
$h_M$ such that $\lim_{t\to 0}M(\eta t)/M(t) = 0$ for some $\eta > 0$. 
Results concerning when an Orlicz space is a Marcinkiewicz space are 
known; see, for instance, \cite{M}. It can be seen that there are 
Orlicz sequence spaces which satisfy the above condition, but are not 
isomorphic to any Marcinkiewicz space.

\begin{lem} \label{dual} 
Let $(e_k)$ be a normalized $1$-symmetric basis of a Banach space $E$, 
and assume 
that $(e_k)$ is not equivalent to the $\ell^1$-basis. Denote by 
$(e'_k)$ the functionals biorthogonal to $(e_k)$. Let $\xi > 0$ be 
given. For any $\ep >
0$, there exists $\delta = \delta(\ep) > 0$ such that if $x \in E$,
\[ \|x\| \geq \xi, \quad \|x\|_\infty \leq \delta, \] 
then there
exists $x' \in \spn\{(e'_k)\}$, 
\[ \|x'\| \leq 1, \quad \|x'\|_\infty \leq
\ep, \quad \text{and} \quad \la x,x'\ra \geq \frac{\xi}{2}.\] 
\end{lem}

\begin{pf} 
Suppose $\ep > 0$ is given.
Since $(e_k)$ is not equivalent to the $\ell^1$-basis, there exists 
$k_0 \in \N$
such that $\|\sum_{k\in A}e'_k\| > 1/\ep$ if $|A| \geq k_0$.  Let
$\delta = \xi(4k_0)^{-1}$.  If $x \in E$, $\|x\| \geq \xi$, and
$\|x\|_\infty \leq \delta$, there exists $\sum b_ke'_k \in 
\spn\{(e'_k)\}$, $\|\sum b_ke'_k\| \leq 1$, such that $\la x, \sum 
b_ke'_k\ra \geq 3\xi/4$.  Let $A =
\{k :  |b_k| > \ep\}$.  Since $\|\sum_{k\in A}e'_k\| \leq
\ep^{-1}\|\sum b_ke'_k\| \leq 1/\ep$, $|A| < k_0$ by the choice of 
$k_0$.
Then 
\[ |\la x, \sum_{k\in A}b_ke'_k\ra| \leq \|x\|_\infty\sum_{k\in
A}|b_k| \leq \delta|A| < \delta k_0 \leq \frac{\xi}{4} .  \] 
If we let $x' =
\sum_{k\notin A}b_ke'_k$, then $x' \in \spn\{(e'_k)\}$, $\|x'\| \leq 
1$, $\|x'\|_\infty \leq
\ep$, and 
\[ \la x, x'\ra \geq \la x, \sum b_ke'_k\ra - |\la x,
\sum_{k\in A}b_ke'_k\ra| \geq \xi/2 .  \] 
This proves the lemma.
\end{pf}

\begin{lem} \label{levels} 
Let $M$ be a nondegenerate
Orlicz function.  Suppose that there exists $\eta > 0$ so that
$\lim_{t\to 0}M(\eta t)/M(t) = 0$.  There is a real sequence
$(\delta_n)^\infty_{n=1}$, strictly decreasing to $0$, such that
whenever $a = (a_k)$ is an element in $h_M$ of norm $1/\eta$, then
either $\|a\|_\infty > \delta_1$, or there exists an $n \in \N$ such 
that
$\|a\chi_{A_n}\| \geq 1$, where $A_n = \{k :  \delta_{n+1} < |a_k|
\leq \delta_n\}$.  
\end{lem}

\begin{pf} 
Choose a strictly decreasing sequence
$(\delta_n)^\infty_{n=1}$ such that $\delta_1 < 1$, and $M(\eta t) <
M(t)/2^{n}$ if $0 < t \leq \delta_n$, $n \in \N$.  If the lemma fails,
there exists $a = (a_k) \in h_M$ of norm $1/\eta$ such that
$\|a\|_\infty \leq \delta_1$ and $\|a\chi_{A_n}\| < 1$ for all $n \in
\N$.  Note that $\supp a \subseteq \cup^\infty_{n=1}A_n$, and 
$\sum_{k\in
A_n}M(|a_k|) < 1$ for all $n$.  Now if $k \in A_n$, then $0 < |a_k| 
\leq
\delta_n$, hence $M(\eta|a_k|) < M(|a_k|)/2^{n}$.  Therefore, 
\begin{align*} 
\sum^\infty_{k=1}M(\eta|a_k|) &=
\sum^\infty_{n=1}\sum_{k\in A_n}M(\eta|a_k|) \\ 
&< \sum^\infty_{n=1}\frac{1}{2^n}\sum_{k\in A_n}M(|a_k|) \\ 
&< \sum^\infty_{n=1}\frac{1}{2^n} = 1 .  
\end{align*} 
But this
contradicts the fact that $\|a\| = 1/\eta$.  The lemma follows.
\end{pf}

\begin{thm} \label{embed} 
Let $M$ be a nondegenerate Orlicz
function.  Suppose that there exists $\eta > 0$ so that $\lim_{t\to
0}M(\eta t)/M(t) = 0$.  Then the Orlicz sequence space $h_M$ is
isomorphic to a subspace of $C(\om^\om)$.  
\end{thm}

\begin{pf} 
As before, let $(e_k)$ be the unit vector basis of $h_M$, and let 
$(e'_k)$ be the sequence of biorthogonal functionals.
The assumption on the Orlicz function $M$ easily implies
that $h_M$ is not isomorphic to $\ell^1$.  Taking $\xi =1$, we obtain 
the function
$\delta$ from Lemma \ref{dual}.  
It follows from the proof of Lemma \ref{levels} that the 
sequence
$(\delta_n)^\infty_{n=1}$ from that lemma can be chosen to satisfy the 
additional condition that
$0 < \delta_n < \delta(2^{-n})$
for all $n \in \N$.  For each $n \in \N$, choose $l_n \in \N$ so that
$\delta_{n+1}\|\sum^{l_n}_{k=1}e_k\| > 1/\eta$.  Let $B_n$ be the
subset of $\R^{l_n}$ consisting of all sequences $(b_k)^{l_n}_{k=1}$
such that 
\begin{gather*} 
\frac{1}{2^n} \geq b_1 \geq b_2 \geq
\dotsb \geq b_{l_n} \geq 0 \\ 
\intertext{and}
\|\sum^{l_n}_{k=1}b_ke'_k\| \leq 1 .  
\end{gather*} 
Finally, let
$\rho_n$ be the norm on $\R^{l_n}$ defined by 
\[\rho_n(a_1,\dotsc,a_{l_n}) = 2^n\,\sup\{\sum^{l_n}_{k=1}a^*_kb_k :
(b_1,\dotsc,b_{l_n}) \in B_n\} .\] 
Clearly, each $\rho_n$ is a
lattice norm.  Moreover, each $\rho_n$ is a normalized norm.
Indeed, since $(1/2^n,0,\dotsc,0) \in B_n$, we see that 
\[
\rho_n(a_1,\dotsc,a_{l_n}) \geq 2^n(\frac{a^*_1}{2^n}) = a^*_1 =
\|(a_1,\dotsc,a_{l_n})\|_\infty .  \] 
On the other hand, 
\[
\rho_n(a_1,\dotsc,a_{l_n}) \leq 2^n\cdot
\frac{1}{2^n}\cdot\sum^{l_n}_{k=1}a^*_k = \|(a_1,\dotsc,a_{l_n})\|_1
,\] 
since $\|(b_1,\dotsc,b_{l_n})\|_\infty \leq 1/2^n$ for all
$(b_1,\dotsc,b_{l_n}) \in B_n$.  This shows that $\rho_n$ is a
normalized lattice norm.  Next, we claim that if $a =
\sum^\infty_{k=1}a_ke_k \in h_M$, then 
\begin{equation} \label{norming} 
\sup_n\frac{1}{2^n}\rho_n(a^*_1,\dotsc,a^*_{l_n}) \leq
\|a\| \leq \max\{\frac{1}{\eta\delta_1}\|a\|_\infty,
\frac{2}{\eta}\sup_n\frac{1}{2^n}\rho_n(a^*_1,\dotsc,a^*_{l_n})\} .
\end{equation} 
Let $n \in \N$.  If $(b_1,\dotsc,b_{l_n}) \in B_n$,
then $\|\sum^{l_n}_{k=1}b_ke'_k\| \leq 1$.  Therefore, 
\[
\rho_n(a^*_1,\dotsc,a^*_{l_n}) \leq 2^n\|\sum^\infty_{k=1}a^*_ke_k\| =
2^n\|a\| .\] 
This proves the left half of (\ref{norming}).  If
$\|a\| = 1/\eta$, define the sets $A_n$ as in Lemma \ref{levels} using
the sequence $(\delta_n)$.  By the same lemma, either $\|a\|_\infty >
\delta_1$, or there exists $n_0$ such that $\|a\chi_{A_{n_0}}\| \geq
1$.  In the latter case, we also observe that
$\|a\chi_{A_{n_0}}\|_\infty \leq \delta_{n_0} < \delta(2^{-n_0})$.
Hence Lemma \ref{dual} yields $x' = \sum^\infty_{k=1}c_ke'_k \in h'_M$
such that 
\[ \|x'\| \leq 1, \quad \|x'\|_\infty \leq
\frac{1}{2^{n_0}}, \quad \text{and} \quad \la a\chi_{A_{n_0}}, x'\ra
\geq \frac{1}{2} .\] 
Notice that since $\delta_{n_0+1}\|\chi_{A_{n_0}}\| \leq
\|a\| = 1/\eta$, $|A_{n_0}| < l_{n_0}$.  Let $c =
(c^*_1,\dotsc,c^*_{l_{n_0}})$.  Since $\|x'\| \leq 1$ and
$\|x'\|_\infty \leq 1/2^{n_0}$, $c \in B_{n_0}$.  But then 
\begin{align*} 
\frac{1}{2^{n_0}}\rho_{n_0}(a^*_1,\dotsc,a^*_{l_{n_0}})
&\geq \sum^{l_{n_0}}_{k=1}a^*_kc^*_k \\ 
&\geq \sum_{k\in A_{n_0}}a_kc_k \\ 
&= \la a\chi_{A_{n_0}}, x'\ra \geq \frac{1}{2}.
\end{align*} 
Thus, if $\|a\| = 1/\eta$, then either $\|a\|_\infty >
\delta_1$ or $\sup_n\frac{1}{2^n}\rho_n(a^*_1,\dotsc,a^*_{l_n}) \geq
1/2$.  The right half of inequality (\ref{norming}) follows
immediately.  From (\ref{norming}), one easily deduces that the norm
on $h_M$ is equivalent to the norm $\rho$ defined by 
\[ \rho(a) =
\sup_n\frac{1}{2^n}\rho_n(a^*_1,\dotsc,a^*_{l_n}) \] 
for $a =
(a_k)$.  We may now apply Theorem \ref{general} to conclude that $h_M$
is isomorphic to a subspace of $C(\om^\om)$.  
\end{pf}

\section{Symmetric sequence spaces}

If $E$ is a Banach space, a subset $W$ of the dual space is {\em
norming}\/ if there is a strictly positive constant $C$ such that 
\[
C^{-1}\sup_{x'\in W}|\la x, x'\ra| \leq \|x\| \leq C\,\sup_{x'\in
W}|\la x, x'\ra| \] 
for every $x \in E$. The proof of the next propostion is left to the 
reader.

\begin{prop} \label{easy} 
A Banach space $E$ is isomorphic to a
subspace of $C(K)$ for some compact Hausdorff topological space $K$ if
and only if $E'$ contains a norming subset $W$ which is a continuous
image of $K$ when endowed with the weak$^*$ topology.  
\end{prop}

To be able to use Theorem \ref{sets}, we would like to shift our focus 
from the norming subset $W$ to the collection of supports of the 
functionals in $W$ (assuming that $E$ is a sequence space). Lemma 
\ref{shorten} shows how to perturb the set $W$ so as to obtain finite 
supports. Then we use Lemma \ref{discrete} to ``discretize'' the 
coordinates of the functionals to ensure that the set of supports is a 
compact subset of $\Ps(\N)$.

\begin{lem} \label{shorten} 
Suppose that $(e_k)$ is a shrinking
unconditional basis of a Banach space $E$, with biorthogonal 
functionals $(e'_k)$, and let $\alpha$ be a
countable ordinal.  Given a continuous function $f : \al +1 \to 
(E',\mbox{weak}^*)$, and $\ep > 0$, there is a continuous map $g :
\al +1 \to (E',\mbox{weak}^*)$ such
that $g(\beta) \in \spn\{(e'_k)\}$ and $\|f(\beta) - g(\beta)\| \leq
\ep$ for all $\beta \leq \alpha$.  
\end{lem}

\begin{pf} 
Induct on $\alpha < \om_1$.  The result is trivial if
$\alpha$ is a finite ordinal.  Now suppose it holds for all $\beta <
\alpha$, where $0 < \alpha < \om_1$.  Assume that  $f : \al +1 \to 
(E',\mbox{weak}^*)$ is continuous, and let $\ep >
0$ be given. First, consider the case when
$\alpha$ is a successor ordinal, say $\alpha = \gamma +1$. By the 
inductive hypothesis, there exists a continuous function 
$g : \gamma +1 \to (E',\mbox{weak}^*)$ such that $g(\beta)\in 
\spn\{(e'_k)\}$, and
$\|f(\beta) - g(\beta)\| \leq \ep$ for all $\beta \leq \gamma$.  Since
$(e_k)$ is shrinking, there exists $y'_\alpha \in \spn\{(e'_k)\}$ such
that $\|f(\alpha) - y'_\alpha\| \leq \ep$.  Extend $g$ by defining 
$g(\al) = y'_\al$. It is clear that $g$ satisfies the requirements of 
the lemma. 

Suppose that $\alpha$ is a limit ordinal.  Since $\alpha < \om_1$,
there is a sequence of ordinals $(\alpha_n)$ which strictly increases
to $\alpha$. Choose $y'_\alpha \in \spn\{(e'_k)\}$
such that $\|f(\alpha) - y'_\alpha\| \leq \ep/2$.  Define $h : \al + 1 
\to (E',\mbox{weak}^*)$ by $h(\beta) =
f(\beta) + y'_\alpha - f(\alpha)$ for all $\beta \leq \alpha$.  $h$ is 
clearly continuous. By the inductive hypothesis, for each $n$, there 
is a continuous function $g_n : \{\beta : \al_{n-1} < \beta \leq 
\al_n\} \to (E',\mbox{weak}^*)$ such that $g_n(\beta) \in
\spn\{(e'_k)\}$, and $\|h(\beta) - g_n(\beta)\| \leq \ep/2^n$
whenever $\alpha_{n-1} < \beta \leq \alpha_n$.  (We interpret
$\alpha_0 < \beta$ to mean $0 \leq \beta$.)  Now define $g : \al +1 
\to (E',\mbox{weak}^*)$ by
\[ g(\beta) = \begin{cases}
g_n(\beta)& \text{if $\al_{n-1} < \beta \leq \al_n, n \in \N$},\\
y'_\al& \text{if $\beta = \al$}.
\end{cases}\]
It is clear that $g(\beta) \in \spn\{(e'_k)\}$ for all $\beta \leq 
\al$. Moreover, $g$ is continuous
at every $\beta < \alpha$.  We claim that it is continuous at
$\alpha$ as well.  This follows from the observations that
$\lim_{\beta\uparrow\alpha}\|g(\beta) - h(\beta)\| = 0$, and
$\lim_{\beta\uparrow\alpha}h(\beta) = h(\alpha) = y'_\alpha$ in the 
weak$^*$ topology.
Finally, if $\beta < \alpha$, 
\begin{align*} 
\|f(\beta) - g(\beta)\|
&\leq \|f(\beta) - h(\beta)\| + \|h(\beta) - g(\beta)\| \\ 
&\leq \|f(\al) - y'_\al\| + \frac{\ep}{2} \leq  \ep, 
\end{align*} 
while $\|f(\al) - g(\al)\| = \|f(\al) - y'_\al\| \leq \ep$ as well.
This completes the
proof of the lemma.  
\end{pf}

\begin{lem} \label{discrete}
Let $A$ be a countable compact subset of $\R$ which contains $0$, and 
suppose $\ep > 0$ is given. Then there is a continuous function $h : A 
\to \R$ which takes only finitely many values, such that $h(0) = 0$, 
and $|h(a) - a| \leq \ep$ for all $a \in A$.
\end{lem}

\begin{pf}
Choose $l$ large enough that $A \subseteq [-l\ep/2,l\ep/2]$. If $-l 
\leq k \leq l$, choose a number $b_k \in [k\ep/2, (k+1)\ep/2) \bs A$. 
For convenience, let $b_{-l-1} = (-l-1)\ep/2$. Define $g : A \to \R$ 
by 
\[ g(a) = b_k \quad \text{if} \quad b_{k-1} < a \leq b_k, \quad -l 
\leq k \leq l.\]
Since $|b_{k} - b_{k-1}| < \ep$, $|g(a) - a| < \ep$ for all $a \in A$.
Suppose $(a_n)$ is a sequence in $A$ which converges to some $a \in 
A$. Choose $k$ such that $b_{k-1} < a \leq b_k$. Since
$a \in A$, but $b_{k} \notin A$, 
$b_{k-1} < a < b_k$. Hence $b_{k-1} < a_n < b_k$ for all large enough 
$n$. This proves that $g$ is continuous. 
Clearly, $g$ takes only finitely many values. Now define the function 
$H$ on the range of $g$ by 
\[ H(g(a)) = \begin{cases}
  0& \text{if $g(a) = g(0)$},\\
  g(a)& \text{otherwise}.
  \end{cases} \]
Finally, let $h = H\circ g$. Clearly, $h$ is continuous on $A$, 
takes only finitely many values, and $h(0) = 0$. If $g(a) = g(0)$, 
then $a$ and $0$ 
lies within the same interval $(b_{k-1},b_k]$, hence $|a| \leq \ep$.
Therefore, $|h(a) - a| = |a| \leq \ep$. If $g(a) \neq g(0)$, then 
$|h(a) - a| = |g(a) - a| \leq \ep$.
\end{pf}

\begin{prop} \label{support} 
Suppose that $(e_k)$ is an
unconditional basis of a Banach space $E$ which embeds into
$C(\om^\al)$ for some countable ordinal $\alpha$.  There exist a
compact, hereditary subset $K$ of $\Ps(\N)$, $K^{(\al+1)} = 
\emptyset$, and a constant $\eta > 0$ 
such that
for every $\sum a_ke_k \in E$, 
\[ \eta\|\sum a_ke_k\| \leq
\sup\bigl\{\|\sum_{k\in A}a_ke_k\| :  A \in K\bigr\}.  \] 
\end{prop}

\begin{pf}
By suitable renorming, we may assume that both the basis $(e_k)$ and 
the sequence of biorthogonal functionals $(e'_k)$ are normalized.
Since $E$ 
embeds into $C(\om^\al)$, the basis $(e_k)$ must be shrinking 
\cite{J}. Using 
Proposition \ref{easy} and Lemma 
\ref{shorten}, we obtain a continuous function $f : \om ^\al + 1 \to 
(E',\mbox{weak}^*)$ such that $f(\beta) \in \spn\{(e'_k)\}$ for all 
$\beta \leq \om^\al$, and that $f(\om^\al + 1)$ is a norming subset of 
$E'$. After appropriate scaling, we may assume that there is a 
constant $\eta > 
0$ such that 
\[ 2\eta\|x\| \leq \sup_{\beta \leq \om^\al}|\la x, f(\beta)\ra| \leq 
\|x\| \]
for all $x \in E$. Each $f(\beta)$ can be expressed in the form 
$f(\beta) = 
\sum_kf_k(\beta)e'_k$. 
For each $k$, denote by $A_k$ the set $\{f_k(\beta) : \beta \leq 
\om^\al\} \cup 
\{0\}$. Then $A_k$ is a countable compact subset of $\R$. By Lemma 
\ref{discrete}, there is a continuous function $h_k : A_k \to \R$ 
which takes only finitely many values, such that $h_k(0) = 0$, and 
$|h_k(a) - a| \leq \eta/3^k$ for all $a \in A_k$. 
Suppose $\beta \leq \om^\al$, then 
\[ \sum|h_k(f_k(\beta)) - f_k(\beta)|\,\|e'_k\| \leq 
\sum^\infty_{k=1}\frac{\eta}{3^k} = \frac{\eta}{2} .\]
It follows that $h(\beta) \equiv \sum h_k(f_k(\beta))e'_k$ converges 
in $E'$, 
and that 
\begin{equation} \label{dist}
\|h(\beta) - f(\beta)\| \leq \eta/2 \quad \text{for all}\quad \beta 
\leq 
\om^\al.
\end{equation}
We claim that $h : \om^\al + 1 \to  (E',\mbox{weak}^*)$
is continuous.
In fact, $h_k\circ f_k : \om^\al + 1 \to \R$ is a continous function 
for every $k$. Since $h(\om^\al+1)$ is bounded by (\ref{dist}), the 
continuity of $h$ follows.
Finally, let 
\[ K = \{A : A \subseteq \supp h(\beta) \quad \text{for some} \quad 
\beta \leq \om^\al\}. \]
Since $h_k(0) = 0$, $\supp h(\beta) \subseteq \supp f(\beta)$ for each 
$\beta \leq \om^\al$. 
Therefore, $K \subseteq \Ps(\N)$. Obviously $K$ is hereditary.\\

\noindent\underline{Claim}. If $\gamma \leq \al$, and $A \in 
K^{(\gamma)}$, then there exists $\beta \in (\om^\al +1)^{(\gamma)}$ 
such that $A \subseteq \supp h(\beta)$.\\
We prove the claim by induction on $\gamma$. If $\gamma = 0$, there is 
nothing to prove. Suppose that $0 < \gamma_0 \leq \al$, and the claim 
holds for all $\gamma < \gamma_0$. Let $A \in K^{(\gamma_0)}$. Assume 
first that $\gamma_0$ is a successor. Then there is a pairwise 
distinct sequence $(A_n)$ in $K^{(\gamma_0 - 1)}$ which converges to 
$A$. By the inductive assumption, there is a sequence $(\beta_n)$ in 
$(\om^\al+1)^{(\gamma_0 - 1)}$ such that $A_n \subseteq \supp 
h(\beta_n)$ for every $n$. Using a subsequence if necessary, we may 
assume that $(\beta_n)$ converges to some $\beta \in \om^\al + 1$. If 
$k \in A$, then $k \in A_n$ for all large enough $n$. Since $\{\la 
e_k,h(\gamma)\ra : \gamma \leq \om^\al\}$ is a finite set, and $\la 
e_k,h(\beta)\ra = \lim_n\la e_k,h(\beta_n)\ra$, $k \in \supp 
h(\beta)$. Hence $A \subseteq \supp h(\beta)$. If $\beta_n = \beta$ 
for infinitely many $n$, then $A_n \subseteq \supp h(\beta_n) = \supp 
h(\beta)$ for infinitely many $n$. But $\supp h(\beta)$ is finite. 
Thus $(A_n)$ has a constant subsequence, contrary to its choice. 
Therefore, $\beta_n \neq \beta$ for all but finitely many $n$. Thus 
\[ \beta \in \bigl((\om^\al+1)^{(\gamma_0-1)}\bigr)^{(1)} = 
(\om^\al+1)^{(\gamma_0)} .\]
Now consider the case when $\gamma_0$ is a limit ordinal. Since $A \in 
K^{(\gamma_0)}$, $A \in K^{(\gamma)}$ for all $\gamma < \gamma_0$. By 
induction, for each $\gamma < \gamma_0$, there exists $\beta_\gamma 
\in (\om^\al+1)^{(\gamma)}$ such that $A \subseteq \supp 
h(\beta_\gamma)$. There is a sequence of ordinals $(\gamma_n)$ which 
strictly increases to $\gamma_0$, such that $(\beta_{\gamma_n})$ 
converges to some $\beta$. It is easy to see that $\beta \in 
(\om^\al+1)^{(\gamma_0)}$. If $k \in A$, then $\la 
e_k,h(\beta_{\gamma_n})\ra \neq 0$. Arguing as before, we see that 
$\la e_k,h(\beta)\ra \neq 0$, i.e., $k \in \supp h(\beta)$. Therefore, 
$A \subseteq \supp h(\beta)$.
This proves the claim.\\

In particular, according to the claim, if $A \in K^{(1)}$, then there 
exists $\beta \in (\om^\al+1)^{(1)}$ such that $A \subseteq \supp 
h(\beta)$. Thus $A \subseteq \supp h(\beta)$, where $\beta \in \om^\al 
+1$. Hence $A \in K$. Therefore $K$ is a closed subset of $\cP(\N)$, 
and must be compact. Using the claim again, we see that any $A \in 
K^{(\al)}$ is a subset of $\supp h(\beta)$ for some $\beta \in 
(\om^\al+1)^{(\al)} = \{\om^\al\}$. So $K^{(\al)}$ is finite, from 
which it follows that $K^{(\al+1)} = \emptyset$. Summarizing, we see 
that $K$ is a compact, hereditary subset of $\Ps(\N)$ such that 
$K^{(\al+1)} = \emptyset$.
Finally, if $\sum a_ke_k \in E$,
\begin{align*}
2\eta\|\sum a_ke_k\| &\leq \sup_{\beta\in \om^\al+1}|\la\sum 
a_ke_k,f(\beta)\ra|. \\
\intertext{Using (\ref{dist}), we see that}
\frac{3}{2}\eta\|\sum a_ke_k\| &\leq \sup_{\beta\in \om^\al+1}|\la\sum 
a_ke_k,h(\beta)\ra| \\
&\leq \sup_{\beta\in \om^\al+1}\|\sum_{k\in\supp h(\beta)}a_ke_k\| \, 
\|h(\beta)\| \\
&\leq \frac{3}{2}\sup_{\beta\in \om^\al+1}\|\sum_{k\in\supp 
h(\beta)}a_ke_k\|,
\end{align*}
since $\|h(\beta)\| \leq \|f(\beta)\| + \eta/2 \leq 1 + 1/2$. 
As $\supp h(\beta) \in K$ for all $\beta \in \om^\al+1$,
\[ \eta\|\sum a_ke_k\| \leq \sup_{A \in K}\|\sum_{k\in A}a_ke_k\| .\]
This proves the proposition.
\end{pf}

\section{Symmetric sequence subspaces of $C(\al)$} \label{subsp}

In this section, we prove the converses to the embedding theorems, 
Theorems \ref{general} and \ref{embed}. In both instances, we rely on 
the 
characterization of the norm provided by Proposition \ref{support}. 
Then 
Theorem \ref{sets} is used to analyze the resulting set $K$. Two norms 
$\|\cdot\|$ and $\rho$ on a vector space $E$ are said to be {\em 
equivalent}\/ if $(E,\|\cdot\|)$ and $(E,\rho)$ are isomorphic via the 
formal identity map.

\begin{thm} \label{converse}
Let $(e_k)$ be a symmetric basis of a Banach space $E$ which embeds 
into $C(\om^\om)$. Then there are a strictly positive real sequence 
$(t_n)$ which decreases to $0$, and for each $n$, a normalized lattice 
norm $\rho_n$ on $\R^n$  such that the norm on $E$ is equivalent to 
\[ \rho(\sum a_ke_k) = \sup_n t_n\rho_n(a^*_1,\dots,a^*_n). \]
\end{thm}

\begin{pf}
Without loss of generality, we may assume that the sequence $(e_k)$ is 
normalized and $1$-symmetric. 
Since $(e_k)$ is unconditional, and $E$ embeds into $C(\om^\om)$, 
$(e_k)$ is shrinking. 
By Proposition \ref{support}, there are a compact hereditary subset 
$K$ of $\Ps(\N)$, $K^{(\om+1)} = \emptyset$, and $\eta > 0$ such that
\[ \eta\|\sum a_ke_k\| \leq
\sup\bigl\{\bigl\|\sum_{k\in A}a_ke_k\bigr\| :  A \in K\bigr\}  \] 
for every $\sum a_ke_k \in E$. By Theorem \ref{sets}, there exist $B 
\in \Pb(\N)$, and $f : \N \to \N$, strictly increasing to $\infty$, 
such that $|A \cap B| \leq f(\min(A\cap B))$ whenever $A \in K$. List 
the elements of $B$ in a strictly increasing sequence $(j_k)$. Then, 
for every $\sum a_ke_k \in E$
\begin{align*}
\eta\|\sum a_ke_k\| & = \eta\|\sum a^*_ke_{j_k}\| \\
&\leq \sup\{\|(\sum a^*_ke_{j_k})\chi_A\| : A \in K\} \\
&= \sup\{\|(\sum a^*_ke_{j_k})\chi_{A\cap B}\| : A \in K\} \\
&\leq \sup_n\bigl\|\sum^{n+f(n)-1}_{k=n}a^*_ke_k\bigr\| ,
\end{align*}
since $|A \cap B| \leq f(\min(A\cap B))$ for all $A \in K$. If 
$\sup_n\|\sum^n_{k=1}e_k\| < \infty$, then $(e_k)$ is equivalent to 
the $c_0$ basis. We can simply take $(t_n)$ to be any strictly 
positive sequence decreasing to $0$, and $\rho_n(a^*_1,\dots,a^*_n) = 
\|(a^*_1,\dots,a^*_n)\|_\infty$ for all $n$.
Otherwise, assume that $\lambda_n \equiv \|\sum^n_{k=1}e_k\| \to 
\infty$ as $n \to \infty$. Using Lemma \ref{dual}, find a 
strictly positive sequence $(t_n)$ which decreases to $0$, such that 
if $x \in E$, $\|x\| \geq \eta/2$, $\|x\|_\infty \leq 1/\lambda_n$, 
then there exists $x' \in E'$, $\|x'\| \leq 1$, $\|x'\|_\infty \leq 
t_n$, and $\la x,x'\ra \geq \eta/4$. Of course, $t_1$ may be chosen to 
be $\leq 1$.
For each $n$, define the lattice norm $\rho_n$ on $\R^{f(n)}$ by
\[ \rho_n(a_1,\dots,a_{f(n)}) = 
\frac{1}{t_n}\sup\bigl\{\sum^{f(n)}_{k=1}a^*_kb_k : \|(b_k)\|_\infty 
\leq t_n,\quad \|\sum b_ke'_k\| \leq 1\bigr\}. \]
It is easy to check that $\rho_n$ is a normalized lattice norm for 
every $n$.
Now if $\|\sum a_ke_k\| = 1$, choose $n$ (depending on $\sum a_ke_k$) 
such that $\|\sum^{n+f(n)-1}_{k=n}a^*_ke_k\| \geq \eta/2$. Note that 
\[ 1 = \|\sum a^*_ke_k\| \geq \|\sum^n_{k=1}a^*_ne_k\| = 
a^*_n\|\sum^n_{k=1}e_k\| . \]
Hence  $a^*_n \leq 1/\lambda_n$. 
By the choice of $t_n$, and the definition of $\rho_n$, 
$\rho_n(a^*_n,\dots,a^*_{n+f(n)-1}) \geq \eta/(4t_n)$. Thus
\[ \rho_n(a^*_1,\dots,a^*_{f(n)}) \geq
\rho_n(a^*_n,\dots,a^*_{n+f(n)-1}) \geq \frac{\eta}{4t_n} .\]
Therefore, for any 
$\sum a_ke_k \in E$,
\[
\frac{\eta}{4}\|\sum a_ke_k\| \leq 
\sup_nt_n\rho_n(a^*_1,\dots,a^*_{f(n)}) .\]
It follows easily that the norm $\|\sum a_ke_k\|$ is equivalent to 
$\sup_nt_n\rho_n(a^*_1,\dots,a^*_{f(n)})$. The conclusion of the 
theorem is now clear.
\end{pf}

We now turn to the converse of Theorem \ref{embed}. Let $M$ be a given 
Orlicz function. If $x = (t_1,t_2,\ldots)$ is a finitely supported 
real sequence, let $\Phi(x) = \sum M(|t_k|)$. Let the countable limit 
ordinal $\al$ be fixed until the end of Proposition \ref{phibound}.
Recall the function $b_\al$, and the sets $\A^f_\beta$, $\beta < \al$, 
associated with $\al$.\\

\noindent {\bf Definition.}\  Suppose that a quadruple $(f,M,S,B)$ is 
given, where $f \in \F$, 
$M$ is a nondegenerate Orlicz function, $S \subseteq (0,\infty)$, 
$\inf S = 0$, and $B \in \Pb(\N)$. Let $a$ be a positive real number. 
A {\em level $0$ block of size $a$} (with respect to $(f,M,S,B)$) is a 
vector of the form $x = t\chi_C$, where $C \in \Ps(B)$, $t \in S$, 
$|C| \geq 1$, and $a/2 \leq \Phi(t) \leq a$. If $\beta < \al$, a {\em 
level $\beta + 1$ block of size $a$} is a vector of the 
form $x = \sump x_i$,, where each $x_i$ is a level $\beta$ block of 
size
$a/p$, $x_1 < x_2 < \dots < x_p$, and $[f(\min(\supp x))]^2 \leq p$. 
Finally, if $\beta \in I_\al$, a {\em level $\beta$ 
block of size $a$} is a vector $x$ such that $x$ is a level 
$b_\al(\beta,f(\min(\supp x)))$ block of size $a$.\\

The reader can easily check that if $x$ is a level $\beta$ block of 
size $a$, then $a/2 \leq \Phi(x) \leq a$. These blocks are constructed 
so that the support of a level $\beta$ block is ``long'' when compared 
with any element in $\A^f_\beta$ which preceeds it. 
Now suppose that a quadruple $(f,M,S,B)$ is
given.  Assume additionally that 
\begin{align} \label{fi}
\sum^\infty_{i=m}\frac{1}{(f(i))^2} &\leq \frac{2}{(f(m))^2}\,
,\quad\text{and} \\  
\frac{f(m)}{(f(m+1))^2} &\leq
\frac{1}{4} \label{fm}
\end{align} 
for all $m \in \N$.

\begin{prop} \label{phibound}
Suppose that  $\beta < \al$, and  $x$ is a level $\beta$ 
block of size $a$, $A \in \A^f_\beta$, and $\min A < \min(\supp x)$. 
Then
\[ \Phi(x\chi_A) \leq \frac{2a f(\min A)}{[f(\min(\supp x))]^2} . \]
\end{prop}

\begin{pf}
The proposition clearly holds for $\beta = 0$, since in this case the 
sets $A$ and $\supp x$ are disjoint.  Assume it holds for some 
$\beta < \al$. We want to show that it also holds for $\beta +1$. 
Suppose $x$ is a level $\beta + 1$ block of size $a$, $A \in 
\A^f_{\beta 
+1}$, and $\min A < \min(\supp x)$. Write $A = \cup^k_{i=1}A_i$, where 
$A_1 < \dots < A_k$, $A_1, \dots, A_k \in \A^f_\beta$, and $k \leq 
f(\min 
A)$. Also $x = \sump x_i$, where each $x_i$ is a level $\beta$ block 
of 
size $a/p$, $x_1 < \dots < x_p$, and $[f(\min(\supp x))]^2 \leq p$. 
For each $A_j$, let 
\[ I_j = \{i : 1 \leq i \leq p, \quad \supp x_i \cap A_j \neq 
\emptyset\}. \]
For each $j$, there is at most one $i_j \in I_j$ such that $\min(\supp 
x_{i_j}) \leq \min A_j$. If no such $i_j$ exists, we let $x_{i_j} = 
0$. Then 
\[ x\chi_{A_j} = x_{i_j}\chi_{A_{j}} + 
\sum\begin{Sb} i\in I_j\\ i \neq i_j \end{Sb}x_i\chi_{A_j} . \]
Notice that if $i \in I_j$, $i \neq i_j$, then $\min A_j < \min(\supp 
x_i)$. Therefore, the inductive assumption applies to 
$x_i$ and $A_j$. We conclude that 
\begin{align*}
\sum\begin{Sb} i\in I_j\\ i \neq i_j \end{Sb}\Phi(x_i\chi_{A_j})
 &\leq \sum\begin{Sb} i\in I_j\\ i \neq i_j \end{Sb}
 \frac{2a f(\min A_j)}{p[f(\min(\supp x_i))]^2} \\
 &\leq \frac{2a f(\min A_j)}{p}\sum^\infty_{i=\min 
A_j+1}\frac{1}{(f(i))^2} \\
 &\leq \frac{4a f(\min A_j)}{p [f(\min A_j +1)]^2} \leq \frac{a}{p} .
\end{align*}
Here we have used the growth conditions (\ref{fi}) and (\ref{fm})
on the last two inequalities. Since $\Phi(x_{i_j}\chi_{A_j}) \leq 
\Phi(x_{i_j}) \leq a/p$, it follows that $\Phi(x\chi_{A_j}) \leq 
2a/p$. Therefore,
\[ \Phi(x\chi_A) = \sum^k_{j=1}\Phi(x\chi_{A_j}) \leq \frac{2ka}{p} 
\leq \frac{2a f(\min A)}{[f(\min(\supp x))]^2} .\]

Finally, suppose that $\beta \in I_\al$, and the 
proposition holds for all ordinals $\gamma < \beta$.. If $A \in 
\A^f_\beta$, 
and $x$ is a level $\beta$ block of size $a$ such that $\min A < 
\min(\supp x)$, then $A \in \A^f_{b_\al(\beta,f(\min A))}$, and $x$ is 
a 
level $b_\al(\beta,f(\min(\supp x)))$ block of size $a$. Let $m = 
f(\min 
A)$, and $n = f(\min(\supp x))$. Then $m < n$. Thus $A \in 
\A^f_{b_\al(\beta,n)}$ by Theorem \ref{order}. Since the proposition 
is 
assumed to hold for the 
ordinal $b_\al(\beta,n)$, we see that
\[ \Phi(x\chi_A) \leq \frac{2a f(\min A)}{[f(\min(\supp x))]^2} \,. \]
This completes the proof of the proposition.
\end{pf}

\begin{prop} \label{main} 
Let $M$ be a nondegenerate Orlicz function.
Suppose that there exist a compact, hereditary subset $K$ of 
$\Ps(\N)$, and a
constant $\eta > 0$ such that 
\begin{equation} \label{nm}
\eta\|a\|_M \leq \sup\{\|a\chi_A\|_M :  A \in K\} 
\end{equation} 
whenever $a \in h_M$.  Then $\lim_{t\to 0}M(\eta t/2)/M(t) = 0$.
\end{prop}

\begin{pf} 
Assume that the proposition fails. Then there are a set $S \subseteq 
(0,\infty)$, $\inf S = 0$, and a $\theta > 0$ such that $M(\eta t/2) 
\geq \theta M(t)$ for all $t \in S$. Since $K$ is countable compact, 
there is a countable ordinal $\beta$ 
such that $K^{(\om^\beta+1)} = \emptyset$. Choose a countable limit 
ordinal $\al$ such that $\beta < \al$. Define the sets $A^f_\gamma$ 
for all $\gamma < \al$ using the function $b_\al$. By Theorem 
\ref{sets}, there 
exist $B \in \Pb(\N)$, and a 
function $f \in \F$, such that $A \cap B 
\in \A^f_\beta$ for all $A \in K$. 
Because of Proposition \ref{fg}, 
it may be assumed that the growth conditions (\ref{fi}) and (\ref{fm}) 
hold. 
Using the quadruple $(f,M,S,B)$ thus 
obtained, construct vectors $x_1 < x_2 < \dots$ such that $x_k$ is a 
level $\beta$ block of size $1$ for every $k \in \N$.
If $A \in K$, then $A \cap B \in \A^f_\beta$. Let $I = \{k \in \N : 
\supp x_k \cap A \neq \emptyset\}$. There is at most one $k_0 \in I$ 
such that $\min(\supp x_{k_0}) \leq \min(A \cap B)$. If $k \in I$, $k 
\neq 
k_0$, then $\min(A\cap B) < \min(\supp x_k)$. Appeal to Proposition 
\ref{phibound} to see that 
\[ \Phi(x_k\chi_A) = \Phi(x_k\chi_{A\cap B}) \leq \frac{2 f(\min(A\cap 
B))}{[f(\min(\supp x_k))]^2} .\]
On the other hand, $\Phi(x_{k_0}\chi_A)  
\leq \Phi(x_{k_0}) \leq 1$. Thus,
\begin{align*}
\sum_{k\in I}\Phi(x_k\chi_A) &=
\Phi(x_{k_0}\chi_A) + \sum\begin{Sb} k \in I \\ k \neq k_0 
\end{Sb}\Phi(x_k\chi_A) \\
&\leq 1 + 2 f(\min(A\cap B))\sum^\infty_{i=\min(A\cap 
B)+1}\frac{1}{(f(i))^2} \\
&\leq 1 + \frac{4 f(\min(A\cap B))}{[f(\min(A\cap B)+1)]^2} \leq 2. 
\end{align*}
Hence, for all $j < \infty$, $\Phi((\frac{1}{2}\sum^j_{k=1}x_k)\chi_A) 
\leq 1$. It follows that $\|(\frac{1}{2}\sum^j_{k=1}x_k)\chi_A\| \leq 
1$. From the hypothesis, we conclude that 
$\frac{\eta}{2}\|\sum^j_{k=1}x_k\| \leq 1$. 
This in turn implies that
\[ \sum^j_{k=1}\Phi(\frac{\eta}{2}x_k) = 
\Phi(\frac{\eta}{2}\sum^j_{k=1}x_k) \leq 1 \]
for all $j$. However, since every nonzero coordinate of $x_k$ belongs 
to $S$, $\Phi(\eta x_k/2) \geq \theta\Phi(x_k) \geq \theta/2$. We have 
reached a contradiction.
\end{pf}

\begin{thm} \label{last}
The following statements are equivalent for every
nondegenerate Orlicz function $M$:  
\begin{enumerate} 
\item There
exists a constant $\eta > 0$ such that $\lim_{t\to 0}M(\eta t)/M(t) =
0$, 
\item The Orlicz sequence space $h_M$ embeds into $C(\om^\om)$,
\item The Orlicz sequence space $h_M$ embeds into $C(\alpha)$ for some
countable ordinal $\alpha$.  
\end{enumerate} 
\end{thm}

\begin{pf} 
The implication (1) $\Rightarrow$ (2) is Theorem
\ref{embed}.  (2) $\Rightarrow$ (3) is obvious.  If (3) holds,
according to Proposition \ref{support}, there are a compact, 
hereditary subset $K$
of $\Ps(\N)$, and an $\eta > 0$ such that inequality (\ref{nm}) of
Proposition \ref{main} holds.  Condition (1) follows from the same
proposition.  
\end{pf}


\end{document}